\documentclass[a4paper,11pt]{amsart}
\usepackage{amsfonts,amsmath,mathrsfs,amssymb,amsthm,amscd,latexsym,amstext,amsxtra}
\usepackage[usenames,dvipsnames]{color}

\usepackage[left=2cm,right=2cm, top=3cm, bottom=3cm]{geometry}

\usepackage[all,cmtip]{xy}
\usepackage{enumerate}
\xyoption{all}
\usepackage{srcltx} 
\usepackage{textcomp}
\usepackage{parskip}

\usepackage{amsthm, amsmath, amssymb,latexsym} 
\usepackage{amsfonts,epsfig,amscd}
\usepackage[]{fontenc}
\usepackage{xy}
\usepackage{enumerate}
\usepackage[]{fontenc}
\usepackage[all]{xy}

\newtheorem{theorem}{Theorem}[section]
\newtheorem{lemma}[theorem]{Lemma}
\newtheorem{corollary}[theorem]{Corollary}
\newtheorem{proposition}[theorem]{Proposition}
\newtheorem{remark}[theorem]{Remark}
\newtheorem{definition}[theorem]{Definition}

\newcommand{\nc}{\newcommand}
\nc{\cH}{{\mathcal H}}
\nc{\cA}{{\mathcal A}}
\nc{\cG}{{\mathcal G}}
\nc{\cC}{{\mathcal C}}
\nc{\cD}{{\mathcal D}}
\nc{\cO}{{\mathcal O}}
\nc{\cI}{{\mathcal I}}
\nc{\cB}{{\mathcal B}}
\nc{\cY}{{\mathcal Y}}
\nc{\cK}{{\mathcal K}}
\nc{\cX}{{\mathcal X}}
\nc{\cS}{{\mathcal S}}
\nc{\cE}{{\mathcal E}}
\nc{\cF}{{\mathcal F}}
\nc{\cZ}{{\mathcal Z}}
\nc{\cQ}{{\mathcal Q}}
\nc{\cN}{{\mathcal N}}
\nc{\cP}{{\mathcal P}}
\nc{\cL}{{\mathcal L}}
\nc{\cM}{{\mathcal M}}
\nc{\cT}{{\mathcal T}}
\nc{\cW}{{\mathcal W}}
\nc{\cU}{{\mathcal U}}
\nc{\cJ}{{\mathcal J}}
\nc{\cV}{{\mathcal V}}
\nc{\bH}{{\mathbb H}}
\nc{\bA}{{\mathbb A}}
\nc{\bG}{{\mathbb G}}
\nc{\bC}{{\mathbb C}}
\nc{\bO}{{\mathbb O}}
\nc{\bI}{{\mathbb I}}
\nc{\bB}{{\mathbb B}}
\nc{\bY}{{\mathbb Y}}
\nc{\bK}{{\mathbb K}}
\nc{\bX}{{\mathbb X}}
\nc{\bS}{{\mathbb S}}
\nc{\bE}{{\mathbb E}}
\nc{\bF}{{\mathbb F}}
\nc{\bZ}{{\mathbb Z}}
\nc{\bQ}{{\mathbb Q}}
\nc{\bN}{{\mathbb N}}
\nc{\bP}{{\mathbb P}}
\nc{\bL}{{\mathbb L}}
\nc{\bM}{{\mathbb M}}
\nc{\bT}{{\mathbb T}}
\nc{\bW}{{\mathbb W}}
\nc{\bU}{{\mathbb U}}
\nc{\bD}{{\mathbb D}}
\nc{\bJ}{{\mathbb J}}
\nc{\bV}{{\mathbb V}}
\nc{\bbZ}{{\mathbb Z}}
\nc{\bR}{{\mathbb R}}
\nc{\fm}{{\mathfrak m}}
\newcommand{\dra}{\dashrightarrow}

\nc{\fr}{{\rightarrow}}
\nc{\co}{{\nabla}}

\nc{\cu}{{\overlineline{\nabla}}}
\nc{\gmc}{\nabla}
\nc{\mtin}[1]{\mbox{{\tiny #1}}}
\nc{\rank}[1]{r_{\mbox{{\tiny #1}}}}

\nc{\OLD}{\textcolor{blue}{OLD. (}}
\nc{\NEW}{\textcolor{blue}{NEW. (}}
\nc{\blue}[1]{\textcolor{blue}{ #1}}

\DeclareMathOperator{\de}{d}
\DeclareMathOperator{\Aut}{Aut}
\DeclareMathOperator{\Bij}{Bij}

\DeclareMathOperator{\Ext}{Ext}
\DeclareMathOperator{\Hom}{Hom}

\DeclareMathOperator{\Ima}{Im}

\DeclareMathOperator{\rk}{rk}

\DeclareMathOperator{\spec}{Spec}

\title{Massey Products and Fujita decompositions \\ on fibrations of curves} 																											
\makeatletter

\@addtoreset{equation}{section}
\makeatother
\author[P. Pirola]{Gian Pietro Pirola}
\address{Dipartimento di Matematica,
	Universit\`a di Pavia,
	Via Ferrata, 1,
	27100 Pavia, Italy}
\email{gianpietro.pirola@unipv.it }

\author[S. Torelli]{Sara Torelli}
\address{Dipartimento di Matematica,
	Universit\`a di Pavia,
	Via Ferrata, 1,
	27100 Pavia, Italy}
\email{sara.torelli7@gmail.com }

\thanks{The authors were supported by PRIN 2015  "Moduli spaces and Lie Theory", INdAM - GNSAGA and by Ministero dell'Istruzione dell'Universit\`a e della Ricerca: Dipartimenti di Eccellenza Program (2018-2022) - Dept. Univ. of Pavia.}

 

\date{\today}

\begin{document}	
	
	\maketitle
	
	\bigskip
		
		\begin{abstract} Let $f:S\to B$ be a fibration of curves and let $f_\ast\omega_{S/B}=\cU\oplus \cA$ be the second Fujita decomposition of $f.$ In this paper we study a kind of Massey products, which are defined as infinitesimal invariants by the cohomology of a curve, in relation to the monodromy of certain subbundles of $\cU$. The main result states that their vanishing on a general fibre of $f$ implies that the monodromy group acts faithfully on a finite set of morphisms and is therefore finite.
			In the last part we apply our result in terms of the normal function induced by the Ceresa cycle. On the one hand, we prove that the monodromy group of the whole $\cU$ of hyperelliptic fibrations is finite (giving another proof of a result due to Luo and Zuo). On the other hand, we show that the normal function is non torsion if the monodromy is infinite (this happens e.g. in the examples shown by Catanese and Dettweiler).
		\end{abstract}
		
	\section*{Introduction}
	

	
	

In this paper we study some geometrical properties of fibrations $f:S\to B$ between a surface $S$ and a curve $B$ that are both smooth, complex and, unless otherwise specified, projective. These are families of projective curves on $B$ whose general fibre $F$ is smooth of geometric genus $g$. We assume that $g\geq 2$ and we will from now on refer to this by saying that $f$ is a fibration of curves of  genus $g\geq 2$ on $B$. 
More precisely, we look at the direct image sheaf $f_\ast\omega_{S/B}$ of the relative dualizing sheaf $\omega_{S/B},$ which is a vector bundle of rank $g$ (see for instance  \cite{Fuj78a}), and combining some infinitesimal and global techniques associated to the geometric variation of the Hodge structure, we relate an infinitesimal invariant introduced in \cite{C-P_TheGriffiths_1995} with a certain subbundle introduced in \cite{Fuj78b} and we obtain some global information on the monodromy of this variation and more deeply on the fibration itself.
 In \cite{Fuj78a} and \cite{Fuj78b}, Fujita introduced two direct sum decompositions of vector bundles $f_\ast\omega_{S/B}=\cO^{\oplus q_f}_B\oplus \cE=\cU\oplus \cA$, where $\cO^{\oplus q_f}_B$ is a trivial bundle of rank the relative irregularity $q_f$ of $f$, $\cE$ is such that $h^1(B,\cE(\omega_B))=0$, $\cU$ is unitary flat and $A$ is ample, which are compatible under the inclusion $\cO^{q_f}\subset \cU$ (see also \cite{kol87}).  We will refer to these as the {\em first} and {\em second Fujita decomposition}, respectively. 
Despite the original motivation, which wondered about semiampleness, the answer recently given by Catanese and Dettweiler in  \cite{CatDet_TheDirectImage_2014} (see also \cite{CD:Answer_2017} and \cite{CatDet_Vector_2016}) showed the existence of fibrations whose unitary flat bundles has infinite monodromy. 
A first motivation to our work has been suggested by these examples, which capture a new piece of information on the variation of the Hodge structure. Indeed, by the theorem of the fixed part proved by Deligne (see \cite{Del_Theorie_1971} ) it follows that when $\cU$ has finite monodromy then it naturally defines a constant Hodge substructure in the geometric variation of the Hodge structure of the fibers, after eventually performing a finite base change. This is no longer true for these examples where the non-finiteness of the monodromy implies that the flat bundle\ \ $\cU\oplus \overline\cU\subset R^1f_\ast\bC\otimes \cO_B$ cannot be defined over the rational field. This fact has opened up to some recent works concerning Hodge loci and Shimura varieties, e.g. \cite{FGP},   \cite{Pet17} and \cite{Mon10}. We also note that in case of fibrations on an elliptic curve the monodromy is always finite (see e.g. \cite{barja-fujita}) so what is said makes sense when the curve $B$ has genus at least $2$.

 We face the study of the monodromy group
  by studying a certain kind of vanishing on a second order cohomological object in the cohomology of a general fibre of a fibration of curves that is related to the Griffith's infinitesimal invariant through a formula given in \cite{C-P_TheGriffiths_1995}.
  The idea that this could force the monodromy to be finite has been suggested by looking to the behaviour of hyperelliptic fibrations where the monodromy is finite (see \cite[Theorem A.1]{LuZuo_OnTheSlope_2017}) and the canonical normal function induced by the Ceresa cycles of the fibres vanishes (see subsection \ref{SubSec-NonVanCriteria}). 

These {\em  adjoint images} introduced in \cite{C-P_TheGriffiths_1995} (see also \cite{Gonz_OnDef_2016}, \cite{P-Z_Variations_2003}, \cite{RizZuc_Generalized_2017}, \cite{R_Infinitesimal_2008} for several applications) has been interpreted  in \cite{R_Infinitesimal_2008} as Massey products of the Dolbeault complex of the cohomology of a curve of genus $g\geq 2$ and we follow this terminology.
Let us shortly recall its construction (that we give in detail in subsection \ref{Sec-MPOnFibrations}). We consider the cohomology of a general fibre $F$ of $f$  and since it is a differential graded algebra then one can define a {\em Massey-product} $\fm_{\xi}(s_1,s_2)$ for a pair $(s_1,s_2)$ of independent elements in the kernel of the cup product map $\cup {\xi}:H^0(\omega_{F})\to H^1(\cO_{F})$ with the corresponding first order deformation $\xi\in H^1(T_{F})$ (Definition \ref{Def-MTpair}).
The definition extends  to families by considering (local) sections of the kernel $\cK_\partial=\ker \partial$ of the connecting morphism $\partial:f_\ast \Omega^1_{S/B}\to R^1f_\ast\cO_F\otimes \omega_B $ defined by pushing forward the short exact sequence  \begin{equation*}
\xymatrix@!R{
	{0}  & {f^\ast\omega_B}  & {\Omega^1_S}  & {\Omega^1_{S/B}}  & {0,}                                  & 
	\ar"1,1";"1,2"\ar"1,2";"1,3"\ar"1,3";"1,4"\ar"1,4";"1,5"
	\hole
}
\end{equation*}
where $\Omega^1_{S/B}$ is the sheaf of relative differentials of $f$. This is because the restriction of the morphism to a general fibre gives exactly the cup product with the Kodaira-Spencer class $\xi$. 
%

We restrict ourself to the study of vanishing for families of Massey-products given by local flat sections of the unitary flat bundle $\cU$ of the second Fujita decomposition of $f$, through the inclusion $\cU\hookrightarrow \cK_{\partial}$.  More precisely, we consider a subbundle $\cM\subset \cU$ that is Massey-trivial generated. By definition (see Definition \ref{Def-MTSubspSec}), $\cM$ is a flat subbundle whose general fibre $M$ is generated under the monodromy action of $\cU$ by a subspace $W$ of local flat sections of $\cU$ with the property of vanishing on Massey-products  (see Definition \ref{Def-MTSubspSec}).

The main results of the paper concern  the monodromy of $\cU$ and more precisely of those subbundles $\cM$  of $\cU$ that are generated by Massey-trivial subspaces. 
We denote with $\cM(C,C')$ the set of morphisms between two smooth projective curves $C$ and $C'$. Roughly speaking, the content of the main result is that the studied property allows to describe the monodromy group through a subset $\mathscr{K}$ of $\cM(F,\Sigma)$, where $\Sigma$ has genus $g_{\mtin{$\Sigma$}} \geq 2$. We recall that by a classical de Franchis' theorem (see e.g.  \cite{Martens_Obervations_1988} and  \cite{AlzatiPirola_Some_1991}) this set is finite. More precisely, we state the following 
\begin{theorem}\label{Thm-MainG} Let $f:S\to B$ be a semistable fibration of curves of genus $g\geq 2$ on a smooth projective curve $B$ and let $f_\ast\omega_{S/B}=\cA\oplus\cU$ be the second Fujita decomposition. Consider a unitary flat subbundle $\cM$ of $\cU$ that is Massey-trivial generated. Then the monodromy group $G_{\mtin{$\cM$}}$ of  $\cM$ acts faithfully on a subset $\mathscr{K}$ of $\cM(F,\Sigma)$, where $\Sigma$ has genus $g_{\mtin{$\Sigma$}} \geq 2$.
\end{theorem}

\begin{theorem}\label{Thm-MainSbis} Under the assumptions of Theorem \ref{Thm-MainG}, then
	\begin{itemize}
		\item[$(i)$] the monodromy group $G_{\mtin{$\cM$}}$ is finite;
		\item[$(ii)$] any fiber  of $\cM$ is isomorphic to $\sum_{k\in \mathscr{K}} k^*H^0(\omega_{\Sigma})\subset H^0(\omega_F).$ 
		\end{itemize}In other words, after a finite base change provided by $G_{\mtin{$\cM$}}$ the pullback of $\cM$ is trivial of fibre $\sum_{k\in \mathscr{K}} k^*H^0(\omega_{\Sigma})$.
\end{theorem}
Releasing the assumption of semistability, we have the following
\begin{corollary}\label{cor:main} Let $f:S\to B$ be a fibration of curves of genus $g\geq 2$ on a smooth projective curve $B$ and let $\cM\subset \cU$ be a Massey-trivial generated subbundle. Then the monodromy group $G_{\mtin{$\cM$}}$ of  $\cM$ is finite.
	\end{corollary} 
The basic examples of Massey trivial bundles are given by trivial fibrations $f:S=F\times B\to B$ such that $f$ is the projection on $B$ and the theorem states that after a base change we get a fibred map on one of these. Restating the situation in terms of abelian varieties, one could say that the finiteness of the monodromy in general allows to define (up to base changes) a fibred map $S\to A\times B$ on the product of $B$ with an abelian variety $A$ that lies in the kernel of the Albanese morphism between the Albanese varieties of $S$ and $B.$ In the case of the theorem, $A$ is defined by the curve $\Sigma$.  

The main tool of our proof is a lifting property (Lemma \ref{Lem-CharFlatSecN}) where we describe the local system $\bU$ underlying $\cU$ as the image of the sheaf $f_\ast\Omega^1_{S,d}$ via a map $f_\ast\Omega^1_{S,d}\to R^1f_\ast\bC$ that is defined by the external differential (\ref{Dia-HolDeRhams}), where $\Omega^1_{S,d}$ is the sheaf of holomorphic closed forms on $S$ and $R^1f_\ast\bC$ is the sheaf whose restriction to the smooth locus of $f$ is the local system associated to the variation of the Hodge structure. Another important tool is given by a revisited version of the classical Castelnuovo-de Franchis theorem (see e.g. \cite{Bea96}).
The relation with the Massey-trivial condition lies in Proposition \ref{Prop-UMTtoIsotropic}. 
This is specific of fibrations of curves but it could be interesting to find an analogous condition in a higher dimension where a Castelnuovo de-Franchis Theorem has been 
proved in \cite{Cat_Moduli_1991}. 

In the last section we apply our results. A first application is given by relating Massey products to the Griffiths infinitesimal invariant of the  canonical normal function through a formula given in \cite{C-P_TheGriffiths_1995}, as already mentioned. The statement is that infinite monodromy implies that the canonical normal function is non torsion.
A second application answers a question of semiampleness getting  a local condition (the vanishing of the Massey product) as a sufficient condition to have such a (global) property. In particular this happens on families of hyperelliptic curves and gives a new proof of \cite[Theorem A.1]{LuZuo_OnTheSlope_2017}.

In the spirit of understanding the geometric information provided by $\cU$ on the geometric variation of the Hodge structure we conclude the introduction by mentioning two interesting facts. The first one proved in  \cite{FGP} states that the moduli image of the fibers is contained in a proper Hodge locus. The second one concerns a study done in  \cite{CFG15} for the the second fundamental form that allows to obtain some information on the totally geodetic locus. It is a work in progress and objective of a future study that of applying the techniques developed in this article to these since our vanishing condition is expected in certain cases depending on the dimension of the loci.

The paper is organized as it follows. In section 
\ref{Sec-LocSyst} we fixes our preliminaries about local systems, flat bundles, monodromy representations on curves, recalling the specific constructions that come from fibrations,  the geometric variation of the Hodge structure and the second Fujita decomposition; in section \ref{Sec-MPOnFibrations} we recall the definition of Massey products in the setting of families of curves and we improve it by proving a splitting for the kernel of the IVHS; in section \ref{Sec-LiftingsOnU} we prove the lifting property describing $\cU$ (lemma \ref{Lem-CharFlatSecN}) and in section \ref{Sec-MTsubbundles} we link such a description to the vanishing condition on Massey products, by using a revisited version of the classical Castelnuovo de Franchis theorem for fibrations of curves; in section \ref{Sec-ProofMainTheorems} we prove our main results; in section \ref{Sec-Applications} we apply our results.

{\bf Acknowledgements.} The authors would like to thank Fabrizio Catanese for the very helpful discussions on the subject, during which he pointed out an inaccuracy in the proof. In the same spirit we thank Alberto Collino and Enrico Schlesinger for their suggestions and comments. Moreover, we would like to thank the organizers of a series of seminaries on the topic, Elisabetta Colombo, Paola Frediani, Alessandro Ghigi, Ernesto Mistretta, Matteo Penegini, Lidia Stoppino. A last thank goes to Victor {G}onz{\'a}lez {A}lonso, for his contribution to the generalization of the Castelnuovo de Franchis theorem. 

{\bf Assumptions and notations.}\label{SusSec-GeneralitiesOnFibrations}
We work over $\bC$. Let $Y$ be a smooth variety  $Y$ of dimension $r$. Then $T_Y$ is the tangent sheaf of $Y,$  $\Omega^1_Y$ is the sheaf of holomorphic $1$-forms and $\Omega^k_Y=\wedge^k\Omega^1_Y$ is  the sheaf of holomorphic $k$-forms, $\omega_Y$ is the canonical sheaf and since $Y$ is smooth,  $\omega_Y=\wedge^r\Omega^1_Y$. Let $(\Omega^\bullet_Y, d)$ be the holomorphic de Rham complex. We set $\Omega^k_{Y,d}=\ker (d: \Omega^k_Y\to \Omega^{k+1}_Y)$ the subsheaf  $\Omega^k_Y$ of the holomorphic $k-$forms on $Y$ that are closed with respect to $d$. Let $f:S\to B$ be a fibration of curves of genus $g$. Then $\omega_{S/B}=\omega_S\otimes f^*\omega_B^{\vee}$ be the relative dualizing sheaf and  $\Omega^1_{S/B}$ be the sheaf of relative differentials defined as the cokernel of the exact sequence $
	0\to f^*\Omega^1_B \to \Omega^1_S \to \Omega^1_{S/B}\to 0$
induced by the morphism $d f:T_S\to f^*T_B$ associated to the differential. 

\section{Local systems, flat vector bundles and monodromy representations on fibrations of curves.}\label{Sec-LocSyst} In this section we do our preparation about local systems, flat vector bundles, monodromy representations, introducing the reader to the content of the following sections (see \cite{Del70}, \cite{V_HodgeTheoryI_2002}, \cite{V_HodgeTheoryII_2003}, \cite{PetSteen_Mixed_2008}  \cite{Kob_Differential_1987} for more details). 
			
			We consider a smooth irreducible curve $B$, a point $b\in B$ and the fundamental group  $\pi_1(B,b)$ of $B$ with base point $B$. Let $V$ be a finite dimensional complex vector space. We recall  that it is equivalent to give 
				\begin{itemize}
					\item[$(\mtin{LS})$] A {\em Local system} $\bV$  of stalk $V$ on $B$, namely  a sheaf  which is locally isomorphic to the constant sheaf of stalk $V$;   
					\item[$(\mtin{FB})$] A {\em flat vector bundle} $(\cV,\nabla)$  of fibre isomorphic to $V$ on $B$, namely a pair given by a locally free sheaf  $\cV$ of $\cO_B$-modules together with a flat connection $\nabla:\cV\to\cV\otimes\Omega^1_B$ (i.e. such that the curvature $\Theta=\nabla^2$ is identically zero).
			      \item[$(\mtin{MS})$] A {\em Monodromy representation} $\rho_{\mtin{V}}$ of $V$ on  $B$, namely a homomorphism
				$$
				\rho_{\mtin{V}}: \pi_1(B,b)\to \Aut(V).
				$$ 
				\end{itemize}
				
				Shortly, these correspondences are given as
				\begin{itemize}
					\item$\bV\leftrightarrow (\cV, \nabla)$:   $\cV:=\bV\otimes_{\bC}\cO_B$ and $\nabla$ the flat connection defined by $\ker \nabla\simeq \bV;$ 
					\item$\bV\leftrightarrow \rho_{\mtin{V}}$: $\rho_V(\gamma)=\alpha\circ \gamma^*\alpha^{-1},$ where $\gamma^*:\bV_b\simeq  \bV_b$ is the isomorphism induced by $\gamma\in \pi_1(B,b)$ and $\alpha:\bV_b  \simeq V$; 
					   $\bV= \widetilde{B}\simeq V / \rho_V(\pi_1(B,b)),$ where $\widetilde{B}$ is the universal covering of $B$. 
				\end{itemize}

			\begin{definition}\label{def:Mon}
				The {\em Monodromy group} of a flat bundle $\cV$ (or equivalently of the underlying local system $\bV$ ) is the image $\Ima \rho_{\mtin{V}}$ of the monodromy representation accociated to $\cV$.
				\end{definition}
				
				We denote by $H_{\mtin{V}}=\ker \rho_{\mtin{V}}$  the kernel of the monodromy representation $\rho_{\mtin{V}}$  and $G_{\mtin{V}}=\pi_1(B,b)/H_{\mtin{V}}$ the quotient group. Then $G_{\mtin{V}}$ is naturally isomorphic to the monodromy group $\Ima \rho_V$ and we will from now on identify them.

			We say that a flat bundle (equivalently a local system or a monodromy representation) is defined over a field $K\subset \bC$ if there  is a flat bundle $\cV_{\mtin{$K$}}$ whose fibre is isomorphic to a finite vector space $V_{\mtin{$K$}}$ on $K$ such that $V=V_{\mtin{$K$}}\otimes_K \bC$.
		
			\begin{remark} \label{garibaldi}
			There is an isomorphism $\Gamma(A,\bV)\to V$ between the spaces of local flat sections of $\bV$ on an open contractible subset $A\subset B$  and $V$, since $\bV$ is trivial on $A$. We will from now on use this identification without specify it. 
					\end{remark}

We recalll some well known properties about local systems and behaviour of their monoromy groups. 
	\begin{proposition}The category of local systems is stable for 
		\begin{itemize}\item[(i)] sum of two local systems. Moreover, the sum has finite monodromy group when both the addends have it; \item[(ii)] base changes under morphisms $u : B'\to B$ of curves. Moreover, the associated monodromy representation factors through $u_\ast:\pi_1(B',b')\to \pi_1(B,b),$ for $b'\in u^{-1}(b)$, and finiteness is preserved by finite morphisms.  \end{itemize}
		\end{proposition}
%

			A {\em unitary} structure on these objects is an hermitian metric on the vector space that is compatible with the given structure.  
				Unitary flat bundles has been classified by Narasimhan and Seshadri (see \cite{NarSes_Stable_1965}) exactly as those that are stable. By definition, these are those with slope (i.e. the number given by the degree over the rank of a vector bundle) that is not decreasing slope on subbundles.
			
		We now consider a local system $\bV$  of  stalk $V$  with its corresponding monodromy representation $\rho_{\mtin{V}}: \pi_1(B,b)\to \Aut(V)$. We take a vector subspace $W\subset V$ and we want to attach a local subsystem of $\bV$ to $W$. We define $$G_{\mtin{V}}\cdot W:=\sum_{g\in G_{\mtin{V}}} g\cdot W,$$ where $g\cdot W:=\rho_{\mtin{V}}(g)(W)$ (shortly, $gW$). We remark that $G_{\mtin{V}}\cdot W$ is the smallest subspace of $V$ containing $W$ and invariant under the action $\rho_{\mtin{V}}$ so that it defines the smallest subrepresentation of $V$ containing $W.$

			\begin{definition}\label{Def-GenLS} Let $\bV$ be a local system on $B$ of stalk $V$ and  $W\subset V$ be a vector subspace of $V$. The local system $\widehat{\bW}$ {\em generated by $W$} is the local subsystem of $\bV$ of stalk $\widehat{W}=G_{\mtin{V}}\cdot W.$ 
				\end{definition}
				As usual, we denote by $\rho_{\mtin{$\widehat{W}$}}$ the monodromy representation of $\widehat{W},$ with $H_{\mtin{$\widehat{W}$}}$ the kernel and with $G_{\mtin{$\widehat{W}$}}$ the quotient. We also denote by $H_{\mtin{W}}$ the subgroup of $H_{\mtin{V}},$ which fixes pointwise $W.$ We remark that $H_{\mtin{$\widehat{W}$}}$ is the normalization of $H_{\mtin{W}}.$
				
				We have the following properties of generated local systems. 
				\begin{proposition}\label{Prop-SubLocSystMon} Let $\bV$ be a local system of stalk $V$ on $B$. 
					
										\begin{itemize}\item[$(i)$] if $W_1\subset W_2\subset V$ and the local system $\widehat{\bW}_2$ generated by $W_2$ has finite monodromy, then the local system $\widehat{\bW}_1$ generated by $W_1$ has finite monodromy; 
											\item[$(ii)$]if $W \subset \Gamma(A, \bV)$ and $u:B'\to B$ is a Galois covering of curves. Then  local subsystem $\widehat{\bW}$ of $\bV$ generated by $W$ is isomorphic to $\sum_{g_i\in I_u}{\widehat{\bW}}_{g_i}$, 
											where ${\widehat{\bW}}_{g_i}$ is the local subsystem generated by $u^*(g_i\cdot W)$ and $I_u\subset \pi_1(B,b)$ is a set of generators of the cokernel of $u_\ast:\pi_1(B',b')\to \pi_1(B,b).$\end{itemize}
				\end{proposition}
				\begin{proof} $(i)$
					Let $H=\ker\rho$ be the kernel of representation of $\bV,$ and $H_i=\ker \rho_i$ be the kernel of the subrepresentations $\rho_i$ of $\widehat{\bW}_i,$ for $i=1,2$. We have an inclusion $ H_2 \vartriangleleft H_1$  of subgroups which gives a surjection $G_2:=\pi_1(B,b)/H_2 \to G_1:=\pi_1(B,b)/H_1$
					on the quotients groups isomorphic to the monodromy groups of $\widehat{\bW}_2$ and $\widehat{\bW}_1,$ respectively. Thus whenever the monodromy of $\widehat{\bW}_2$ is finite, the monodromy of $\widehat{\bW}_1$ is finite.
					
					$(ii)$ We consider the  inverse image $u^{-1}\widehat{\bW}$ of the local system $\widehat{\bW}$ generated by $W,$ which is a local system of the same stalk (i.e. $\pi_1(B,b)\cdot W$) and monodromy representation $\rho^{-1}_{\mtin{W}}$ given by the action of $\pi_1(B',b')$ through the composition $\rho\circ u_\ast,$  where $u(b')=b$ and  $u_\ast:\pi_1(B',b')\to \pi_1(B,b)$  induced by $u.$ Let $I_u$ be a set of generators for the cokernel of $u_\ast$ and be $\widehat{\bW}_{g}$ be the local system on $B'$ of stalk generated by $u^*gW$ (i.e. $\pi_1(B',b')\cdot u^*gW).$ Since the monodromy action of $g\in \pi_1(B,b)$ sends $W$ to $gW,$ then the sum of these over $I_u$  reconstructs exactly $u^{-1}\widehat{\bW}.$
				\end{proof}

					
			\subsection{Local systems of fibrations of curves: the geometric variation of the Hodge structure and the second Fujita decomposition}\label{SubSec-Prel-LocSystOnFibr} In this section we shortly recall the construction of two local systems, the geometric VHS of weight one and the the unitary flat bundle of  the second Fujita decomposition $f$, which are both defined by the cohomology of a general fibre $F$ of a fbration $f:S\to B$ of curves of genus $g$ (see \cite{Grif_PeriodsIII_1970}, \cite{V_HodgeTheoryII_2003}, \cite{CatElZFouGrif_Hodge_2014} \cite{PetSteen_Mixed_2008} and \cite{Fuj78b}, \cite{CatDet_TheDirectImage_2014}, \cite{CD:Answer_2017}, \cite{CatDet_Vector_2016} \cite{barja-fujita}, respectively, for details). 
			
			{\bf $(1)$ The geometric variation of the Hodge structure of $f$.}\label{ES-LocSyst1}  Let $f:S\to B$ be a fibration of curves of genus $g$ and assume that $f$ is smooth (i.e. $f$ is a submersive morphism and consequently all the fibres are smooth curves). Then $f$ defines the pair $(\bH_{\bZ}=R^1f_*\bZ,\,  \cF^1=f_*\omega_{S/B})$ where
			\begin{itemize}
				\item $\bH_{\bZ}$ is a local system of lattices over $\bZ$;
				\item $\cF^1=f_\ast \omega_{S/B}\subset \cH= R^1f_*\bZ\otimes_{\bC}\cO_B$ is a filtration.
				\end{itemize}			
				This is called geometric variation of the Hodge structure of $f$.
			
			On a fibre $F$ over a point $b\in B$ we have indeed the corresponding weight one Hodge structure  $(H_{\bZ}=H^1(F,\bZ),\,H^{1,0}=H^0(\omega_{F}))$ given by the cohomology of the fibre $F$, where we take the fibres $H^1(F,\bZ)$ of the local system $R^1f_\ast\bZ $ and $H^0(\omega_F)$ of the vector bundle $f_\ast\omega_{S/B}$ that are compatible with the inclusion$ f_\ast\omega_{S/B}\subset R^1f_\ast\bZ \otimes \cO_B$.
			
			Letting $f:S\to B$ acquire some kind of singularities, then the sheaf $R^1f_\ast\bZ$ is not a local system in general depending on the  homology of the singular fibres that can differ from that of a general (smooth) fibre. Anyway, we can always restrict $f$ over the locus $B^0\subset B$ parametrizing the smooth fibres of $f$, which is a Zariski open subset of $B$. Through the inclusion $j:B^0\subset B$ then $f$ defines a VHS as before over $B^0$ setting $(\bH_{\bZ}=j^\ast R^1f_\ast \bZ, j^\ast f_\ast\omega_{S/B})$ together with a morphism 
			\begin{equation}\label{Mor-AdjRestr}
			\alpha : R^1f_\ast\bC_S\to j_\ast j^*R^1f_\ast\bC_S. 
			\end{equation} 
		locally given by restriction to $j^*R^1f_\ast\bC$. It is natural to ask whenever the morphism is an isomorphism and in this case  we say that $R^1f_\ast\bC$ is completely determined by the local system $j^\ast R^1f_\ast \bC.$

Putting together  \cite[Lemma C.13, pag. 440]{PetSteen_Mixed_2008} and \cite[Theorem 5.3.4, pag. 266]{CatElZFouGrif_Hodge_2014}) we get the following 
			\begin{lemma} \label{Lem-LocInvIso}
				Let $f:S\to B$ be a fibration of curves of genus $g$ with isolated singularities (e.g. semistable). Then the morphism \eqref{Mor-AdjRestr} is an isomorphism.	\end{lemma}

			{\bf $(2)$ The second Fujita decomposition of $f.$} Let $f:S\to B$ be a fibration of curves of genus $g$. The {\em second Fujita decomposition} (\cite{Fuj78b},    \cite{CatDet_TheDirectImage_2014}) states  that there is a splitting 
			\begin{equation}\label{Dec-IIFuj}
			f_\ast\omega_{S/B}=\cU\oplus \cA
			\end{equation}
			of $f_\ast\omega_{S/B}$ into the direct sum of a unitary flat bundle $\cU$ and an ample bundle $\cA$. We shortly recall the definition of $\cU$ and its relation with the geometric VHS introduced above.
			According to Section \ref{Sec-LocSyst}, we consider the correspondig unitary local system $\bU$ that has stalk isomorphic to the fibre $U$ of $\cU$ and is such that $\cU=\bU\otimes \cO_B,$ and a unitary monodromy representation $\rho_{U}:\pi_1(B,b)\to \Aut(U,h)$. We denote by $\rho$ (instead of $\rho_{U}$) the monodromy representation, by $H$ the kernel and by $G$ the quotient $\pi_1(B,b)/H.$ We recall that $G$ is naturally isomorphic to the monodromy group of $\cU$ and we identify them. 

			 Through the inclusion $j:B^0\subset B$, we consider the restriction $\cU_{|B^0}=j^\ast \cU$ of 
$\cU$ to $B^0$. This is the unitary flat subbundle of $j^\ast R^1f_\ast\bC\otimes \cO_{B^0}$ defined by taking the larger subspace $U\subset H^0(\omega_F)$ that gives a subrepresentation of $R^1f_\ast\bC$ together with an hermitian form $h$ induced by the topological intersection form $Q$ of the fibers.  So it is naturally defined by the VHS of $f$. The question is if and how the unitary flat structure extends to the whole $B$, namely the relation between $\cU_{|B^0}$ and $\cU$. While the whole geometric variation does not extend on $B$, the unitary structure defined by $h$ obliges the local monodromies defined by the homotopy of the singular fibres to act trivially on $\cU_{|B^0}$ when singularities are simple normal crossing (i.e. $f$ is semistable) (see \cite{CD:Answer_2017}). In other words, in this case $\cU$ is simply the (trivial) extension of $\cU_{|B^0}$. 
%
 In case of worse singularities, the proof of the existence of $\cU$ requires more sofisticated techniques (see \cite{CD:Answer_2017}) but we can apply the semistable reduction theorem to get a base change $u:B'\to B$ given by a finite morphism of curves  ramified exactly over the worse singularities and after eventually taking a resolution of the fiber product 
			\begin{equation}\label{Dia-BaseChangeSemistab}
			\xymatrix@!R{
				{S':=\widetilde{S\times_BB'}} \,        &     {S}              &    \\
				{B'}            \,        & {B}         &    
				\ar"1,1";"2,1"^{f'}   \ar "1,1";"1,2"^>>>>>>{\varphi}
				\ar"1,2";"2,2"^{f} \ar "2,1";"2,2"^{u}
				\hole,
			}
			\end{equation}     
			we produce a semistable fibration $f':S'\to B'$ of curves of genus $g$ and we call $f':S'\to B'$ the {\em semistable-reduceed fibration of $f$}. This fibration of course has its own Fujita decomposition whose unitary flat summand $\cU'$ is described as above. The relation between the unitary summand $\cU$ of $f$ and the unitary summand $\cU'$ of its semistable reduction is the content of the following lemma.
			\begin{lemma}\label{Lem-UunderBC} There exists a short exact sequence 
				\begin{equation}\label{SES-UnitaryFacBaseChange}
				\xymatrix@!R{
				{0} &	{\cK_{\mtin{U}}}  & {\cU'}  & {u^*\cU}  &  {0,} 
					\ar"1,1";"1,2"\ar"1,2";"1,3"\ar"1,3";"1,4"\ar"1,4";"1,5"
					\hole
				}
				\end{equation}
				which is split. Moreover, $\cK_{\mtin{U}}$ is unitary flat and the splitting is compatible with the underlying local systems. 				
				\end{lemma} 
				\begin{proof}
					 We consider the short exact sequence (\cite[Proposition 2.9]{CatDet_TheDirectImage_2014})
				\begin{equation}\label{SES-HodgeBundlesBaseChange}
				\xymatrix@!R{
					{0}  & {f'_\ast\omega_{S'/B'}}  & {u^*f_\ast\omega_{S/B}}  & {\cG}  & {0,} 
					\ar"1,1";"1,2"\ar"1,2";"1,3"\ar"1,3";"1,4"\ar"1,4";"1,5"
					\hole
				}
				\end{equation}
				where $\cG$ is a skyscraper sheaf supported on points over the singular fibers of $f.$ Comparing the second Fujita decompositions of $f$ and $f'$ we get the s.e.s.
					\begin{equation}\label{SES-HodgeBundlesBaseChange}
				\xymatrix@!R{
					{0}  & {\cA'\oplus \cU'}  & {u^*\cA\oplus u^*\cU}  & {\cG}  & {0,}  
					\ar"1,1";"1,2"\ar"1,2";"1,3"^{i'}\ar"1,3";"1,4"\ar"1,4";"1,5"
					\hole
				}
				\end{equation}
				defining a morphism $\cU'\to u^*\cU$ as the composition of the injection $\cA'\oplus \cU' \hookrightarrow u^*\cA\oplus u^*\cU$ with the projection over $u^*\cU$. We prove that the morphism is surjective. First let us notice that the map $\cA'\to u^\ast \cU$ is identically zero since $\cA$ is ample and $u^*\cU$ is unitary flat. Then we use the characterization of unitary flat bundle over curves of genus greater that $2$ (recalled in section \ref{Sec-LocSyst}) to conclude that the quotient is also unitary flat. Since it also must be supported in a skyscreaper sheaf we conclude that it must be zero. The cases of genus $0$ and $1$ are trivial. We refer to \cite{CatDet_TheDirectImage_2014}, \cite{CD:Answer_2017} and also \cite{CatDet_Vector_2016} for details.
					\end{proof}
	By the previous lemma, we get an inclusion $u^\ast \cU\subset \cU'$ that allows to reduce to the semistable case. We will use it during the proof of one of our main theorem.

        \section{Families of Massey products on fibrations of curves }\label{Sec-MPOnFibrations} In this section we recall the construction of the  "Massey-products" for (local) families of curves. We refer to \cite{C-P_TheGriffiths_1995}, \cite{Gonz_OnDef_2016}, \cite{P-Z_Variations_2003}, \cite{NarPirZuc_Poly_2004}, \cite{R_Infinitesimal_2008} for details on Massey products and to  \cite{Grif_InfinitesimalVariationsIII_1983},\cite{GriffithsHarris_Infinitesimal||_1983},\cite{GriffithsHarris_Infinitesimal||i_1983},\cite{V_HodgeTheoryI_2002}, \cite{V_HodgeTheoryII_2003} for details on deformation theory and variation of the Hodge structure.
        
                Let $f:S\to B$ be a fibration of curves of genus $g$ and assume that  $g\geq 2.$ We consider the corresponding first order deformation           \begin{equation}\label{SeS-DeRhamOnSFibre}
         \xymatrix@!R{
         	{\xi:}&{0}  & {\cO_F}  & {\Omega^1_{S|F}}  & {\omega_{F}}  & {0}                                  & 
         	\ar"1,2";"1,3"\ar"1,3";"1,4"\ar"1,4";"1,5"\ar"1,5";"1,6"
         	\hole
         }
         \end{equation} 
         classified by the estension class $\xi\in \Ext^1_{\cO_F}(\omega_F,\cO_F)\simeq H^1(F,T_F)$, which corresponds to the Kodaira-Spencer class.
         The connecting homomorphism $\delta$ associated to the long exact sequence in cohomology
         \begin{equation}\label{LeS-DeRhamOnSFibre}
         \xymatrix@!R{
         	{0}  & {H^0(\cO_F)}  & {H^0(\Omega^1_{S|F})}  & {H^0(\omega_{F})}  & {H^1(\cO_F)},                                  & 
         	\ar"1,1";"1,2"\ar"1,2";"1,3"\ar"1,3";"1,4"\ar"1,4";"1,5"^{\delta=\cup\xi}
         	\hole
         }
         \end{equation}
        is given by the cup product $\cup \xi : H^0(\omega_F)\to H^1(\cO_F).$ By the Griffiths trasversality theorem (see \cite{GrifTopics1984}), this is also the IVHS (in this case induced by the geometric VHS). Let $K_\xi=\ker(\cup\xi)$ be the kernel of the cup procuct.  We consider the map  
         \begin{equation}\label{Mor-Mp/Aj}
         \wedge_\xi\colon \xymatrix@!R{
         	{\bigwedge^2H^0(\Omega^1_{S|F} )}  & {H^0(\bigwedge^2\Omega^1_{S|F})\simeq H^0(\omega_{F})}                  
         	\ar"1,1";"1,2"
         	\hole
         }
         \end{equation}
         defined by composing the wedge product with the isomorphism $\bigwedge^2H^0(\Omega^1_{S|F} )\simeq \omega_F$ induced by sequence \eqref{SeS-DeRhamOnSFibre}. This is called adjoint map (see \cite{C-P_TheGriffiths_1995}).
         
         Lifting a pair $(s_1,s_2)$ of linearly independent elements of $K_{\xi}$ to $H^0(\Omega^1_{S|F})$ we can compute the adjoint map and  take the image, which we denote with $m(s_1,s_2)$. In general the image depends on the chosen liftings but this in well defined and unique in the quotient $H^0(\omega_C)/< s_1,s_2>_\bC$ with the $\bC$-vector space $< s_1,s_2>_\bC$ generated by $(s_1,s_2),$  (since each lifting must differ from another one for an element in $H^0(\cO_F)\simeq\bC$ by \eqref{LeS-DeRhamOnSFibre}).  
         \begin{definition}\label{Def-MPAJ}
         	The equivalence class 
         	\begin{equation}\label{Mor-Mp/Aj}
         	\fm_{\xi}(s_1,s_2):=[m(s_1\,{s}_2)]\in H^0(\omega_F)/< s_1,s_2>_\bC 
         	\end{equation}
         	is called {\em Massey product of $(s_1,s_2)$} along $\xi$. 
         \end{definition}
         
         The natural {\em vanishing} in this setting is the following.
         \begin{definition}\label{Def-MTpair}
         	We say that a pair $(s_1,s_2)\subset K_{\xi}$ is {\em Massey-trivial} (or equivalently {\em has vanishing Massey-tproducts}) if $\fm_\xi(s_1,s_2)=0$ (namely,  $m(s_1,s_2)\in  < s_1,s_2>_\bC$).
         \end{definition}

%
         %
         We now put in family the above situation. We consider the exact sequence 
         \begin{equation}\label{SeS-DeRhamOnFib}
         \xymatrix@!R{
         	{0}  & {f^*\omega_B}  & {\Omega^1_S}  & {\Omega^1_{S/B}}  & {0}                                  & 
         	\ar"1,1";"1,2"\ar"1,2";"1,3"\ar"1,3";"1,4"\ar"1,4";"1,5"
         	\hole
         }
         \end{equation}that is defined by the morphism $d f:T_X\to f^*T_B$, which is induced by the differential of $f$. 
          The quotient $\Omega^1_{S/B}$ is called the sheaf of relative differentials and it is non torsion free in general. By pushing forward it, we get the exact sequence
         \begin{equation}\label{SeS-DeRhamOnFibf_*}
         \xymatrix@!R{
         	{0}  & {f_*f^*\omega_B\simeq\omega_B}  & {f_*\Omega^1_S}  & {f_*\Omega^1_{S/B}}  & {(R^1f_*\cO_S)\otimes\omega_B}  , 
         	\ar"1,1";"1,2"\ar"1,2";"1,3"\ar"1,3";"1,4"\ar"1,4";"1,5"^-{\partial}
         	\hole
         }
         \end{equation}
         whose connecting morphism $\partial$ describes pointwise the IVHS (see \cite{Gonz_OnDef_2016}). We denote by $\cK_\partial=\ker \partial$ the kernel of $\partial.$ As a subsheaf of $f_\ast \Omega^1_{S/B}$ there is no reason in general for it to be torsion free. This happens when the fibration has isolated singularities and in this case it is automatically locally free  since defined on a curve. The fibre of $\cK_{\partial}$ on a regular point $b$ is exactly the kernel of the cup product with the Kodaira-Spencer class of the corresponding (smooth) fibre $F_b.$ 
         Moreover, as a subsheaf of $f_*\Omega^1_{S/B},$ the sheaf $\cK_\partial$ injects in $f_*\omega_{S/B}$ outside the singular locus of $f.$ The relation between the direct images $f_\ast\omega_{S/B}$, which is a vector bundle of rank $g$, and $f_*\Omega^1_{S/B}$ is given by the following.
         
         \begin{proposition}
         	
         	Let $f:S\to B$ be a fibration of curves of genurs $g$ on a smooth curve $B.$ Assume that the all the fibres of $f$ are reduced.  Then $f_*\Omega^1_{S/B}$ is a sheaf of rank $g$ that injects into $f_\ast\omega_{S/B}$ as a sheaf.
         	\end{proposition}
         	\begin{proof} We consider the exact sequence (see \cite{Gonz_PhdTs_2013} or \cite{Gonz_OnDef_2016})
         		\begin{equation}\label{SeS-DirectImages1}
         		\xymatrix@!R{
         			{0}  & {{f^*\omega_B(Z_d)}_{|Z_d}}& {\Omega^1_{S/B}}  & {\omega_{S/B}}  & {{\omega_{S/B}}_{|Z_0}}  & {0,}                                  & 
         			\ar"1,1";"1,2"\ar"1,2";"1,3"\ar"1,3";"1,4"^{\nu}\ar"1,4";"1,5"\ar"1,4";"1,5"\ar"1,5";"1,6"
         			\hole
         		}
         		\end{equation}
         		where $Z=Z_d+Z_0$ is the singular locus of $f,$ with $Z_d$ a divisor and $Z_0$ supported on isolated points. This induces an injective morphism of sheaf  $\nu':(f_*\Omega^1_{S/B})^{\vee\vee}\hookrightarrow f_*\omega_{S/B}$
         		from the double dual $(f_*\Omega^1_{S/B})^{\vee\vee}$ of $f_*\Omega^1_{S/B}$ to $f_\ast \omega_{S/B}$. Then to conclude it is enough to observe that $f_*\Omega^1_{S/B}$ is isomorphic to its double dual when the fibres are reduced.
         		\end{proof}
         So by the previous proposition under the assumption of isolated singularities the sheaf $\cK_\partial$ is also locally free and injects $\nu:\cK_\partial\hookrightarrow f_*\omega_{S/B}$ by restriction. 
         
         We introduce local families of Massey-products. Let $A\subset B^0$ be an open contractible subset of the locus $B^0\subset B$ parametrizing the smooth fibres of $f$.  A first possibility is to work by hand fibrewise. We fix a local trivialization $\sigma\in \Gamma(A,T_B)$ of $T_B$ over $A$ (up to shrinking $A,$ if necessarily) that defines fibrewise by restriction a generator $\sigma_{b}$ of $T_{\Delta_\epsilon,b}^\vee$ together with an isomorphism $\sigma_{b}: \cO_{F_{b}}\otimes  T^{\vee}_{\Delta_\epsilon,b}\simeq \cO_{F_{b}},$ where $\Delta_\epsilon= \spec \bC[\epsilon]/(\epsilon^2)$ is the the ring of dual numbers of the infinitesimal deformation induced on $F_{b}$ by $f.$  This leads fibrewise to the situation of sequence \eqref{SeS-DeRhamOnFib} and thus
%
we can pointwise repeat that construction of Massey-products for a pair of local section $(s_1,s_2)$ in $\cK_\partial$. The procedure defines a section $\fm_{\sigma}(s_1,s_2)\in \Gamma(A, f_*\omega_{S/B})$ that is  well defined modulo the $\cO_B(A)-$submodule $<s_1,s_2>_{\cO_B(A)}$ of $\Gamma(A,f_*\omega_{S/B})$ generated by $s_1,s_2$.
         \begin{definition}\label{Def-MPAJLoc}
         	We say that {\em a local family of Massey-products} of the pair $(s_1,s_2)$ of sections in $\cK_\partial$ along $\sigma$ in $T_B$ is a section 
         	\begin{equation}\label{Mor-Mp/AjLoc}
         	\fm_{\sigma}(s_1,s_2)\in \Gamma(A, f_*\omega_{S/B}) \end{equation}
          defined modulo the $\cO_B(A)-$submodule $<s_1,s_2>_{\cO_B(A)}$ of $\Gamma(A,f_*\omega_{S/B})$. Such a pair is Massey-trivial (or equivalently, trivial along $\sigma$) if it is pointwise Massey-trivial (Definition \ref{Def-MTpair}).
          
         \end{definition}
         	
%

	More intrinsically, we consider the short exact sequence
	\begin{equation}\label{SeS-Kernel}
	\xymatrix@!R{
		{\zeta:} &{0}  & {\omega_B}  & {f_*\Omega^1_S}  & {\cK_\partial}  & {0} 
		\ar"1,2";"1,3"\ar"1,3";"1,4"\ar"1,4";"1,5"\ar"1,5";"1,6"
		\hole
	}
	\end{equation}
	of locally free sheaves of $\cO_B-$modules defined by \eqref{SeS-DeRhamOnFibf_*} and set $\zeta\in \Ext^1_{\cO_B}(\cK_\partial,\omega_B)$ the extension class. We now prove that one can define a splitting on \ref{SeS-Kernel} and then use it to compute Massey-products in family. 
	
         	\begin{lemma}[Splitting on $\cK_{\partial}$]\label{Lem-SesSPLITKER}
         		Let $f:S\to B$ be a semistable fibration on a smooth projective curve $B.$ Then the short exact sequence 
         		\begin{equation}\label{SeS-KernelSplit<}
         		\xymatrix@!R{
         			{\zeta:} & {0}  & {\omega_B}  & {f_*\Omega^1_S}  & {\cK_\partial}  & {0}  , 
         			\ar"1,2";"1,3"\ar"1,3";"1,4"\ar"1,4";"1,5"\ar"1,5";"1,6"\ar@/^-1.0pc/"1,5";"1,4"_<<<\eta
         			\hole
         		}
         		\end{equation}
         		is split, namely $f_*\Omega^1_S\simeq\omega_B\oplus \cK_\partial$.
         	\end{lemma}
         	\begin{proof} 
Through the chain of isomorphisms 
         		$\Ext^1(K_\partial, \omega_B)\simeq \Ext^1(\cO_B,K_\partial^{\vee}\otimes\omega_B)\simeq H^1(K_\partial^\vee\otimes\omega_B)$
         		we look at $\zeta$ of sequence $\eqref{SeS-KernelSplit<}$ as an element in $H^1(K_\partial^\vee\otimes\omega_B).$ We prove that $\zeta$ is zero, which means by definition that the short exact sequence \ref{SeS-KernelSplit<} is split.
         		
         		To do this we consider the long exact sequence in cohomology
         		\begin{equation}\label{Dia-LongExSec1}
         		\xymatrix@!R{
         			{0}  & {H^0(\omega_B)}  & {H^0(f_*\Omega^1_S)}  & {H^0(\cK_\partial)}  & {H^1(\omega_B)}      & {H^1(f_*\Omega^1_S)}  & {} &
         			\ar"1,1";"1,2"\ar"1,2";"1,3"\ar"1,3";"1,4"\ar"1,4";"1,5"^{\delta}\ar"1,5";"1,6"
         			\hole
         		}
         		\end{equation}
         		and the dual map $\delta^\vee: H^0(\cK_\partial)^\vee\to H^1(\omega_B)^\vee$ of $\delta.$ Through the isomorphisms $H^1(\omega_B)^\vee\simeq H^0(\cO_B)$ and $H^0(\cK_\partial)^\vee\simeq H^1(\omega_B\otimes \cK^\vee),$ we rewrite it as the  map $\delta^\vee: H^0(\cO_B)\to H^1(\omega_B\otimes \cK^\vee)
         		$
         		sending sends $1\mapsto \zeta.$ So $\zeta$ is zero if $\delta^\vee$ is identically zero. It is now enough to prove that the dual map $\delta: H^0(\cK_\partial)\to H^1(\omega_B)$ is identically zero to conclude the same for $\delta^\vee$. By the long exact sequence  \ref{Dia-LongExSec1}, this is equivalent to prove the map  $H^1(B,\omega_B)\to H^1(B,f_*\Omega^1_S)$ is an injection.
         		
         		First of all,  the map $H^1(B,\omega_B)\to H^1(B,f_\ast\Omega^1_S)$ is non zero. Indeed the pullback map $H^1(B,\omega_B)\to H^1(S,\Omega^1_S)$ is an injection, sending the class of a point $b$ on $B$ (which corresponds to a K\"{a}hler form) to the class of the fibre $F$ in $S$ (which is non zero). Then the proof follows since the map must factorize through the Leray spectral sequence (see \cite{V_HodgeTheoryII_2003})
         		\begin{equation}\label{Dia-PullBackLeray}
         		\xymatrix@!R{
         			& &  {H^1(B,\omega_B)}    &  & \\
         			{0}                &  {H^1(B,f_*\Omega^1_S)}    & {H^1(S,\Omega^1_S)}  &    {H^0(B,R^1f_*\Omega^1_S).}      
         			\ar"2,1";"2,2"\ar"2,2";"2,3"\ar"2,3";"2,4"
         			\ar@{^{(}->}"1,3";"2,3"   
         			\ar@{^{(}->}"1,3";"2,2" 
         			\hole
         		}
         		\end{equation}
         	\end{proof}

We introduce the sheaf map that locally defines the adjoint map above (\cite{C-P_TheGriffiths_1995}). Consider the morphism
\begin{equation}\label{Mor-MPRelativeForms}
\xymatrix@!R{
	{\phi : \bigwedge^2f_*\Omega^1_S\otimes T_B}                &  {f_*\omega_{S/B},}    
	\ar"1,1";"1,2"
	\hole
}
\end{equation} 
defined by the morphism  $\bigwedge^2f_*\Omega^1_S\to f_*\bigwedge^2\Omega^1_S$ twisted by $T_B.$ By the projection formula, we have  $f_*\bigwedge^2\Omega^1_S\otimes T_B\simeq f_*\omega_{S/B}.$
By looking at the short exact sequence
\begin{equation}\label{SeS-wedge2KernelSPlit}
\xymatrix@!R{
	{0} & {\omega_B\otimes \cK_\partial}  & {\bigwedge^2f_*\Omega^1_S}  & {\bigwedge^2\cK_\partial}  & {0}  
	\ar"1,1";"1,2"\ar"1,2";"1,3"\ar"1,3";"1,4"\ar"1,4";"1,5"
	\hole
}
\end{equation}
induced by $\zeta,$ which is split by lemma \ref{SeS-KernelSplit<}, we get an injection $\bigwedge^2\cK_\partial \hookrightarrow\bigwedge^2f_*\Omega^1_S$
that factorizes through the morphism \ref{Mor-MPRelativeForms}. Thus we obtain
\begin{equation}\label{Mor-MPRelativeFormsKernel}
\xymatrix@!R{
	{\phi : \bigwedge^2\cK_\partial\otimes T_B}                &  {f_*\omega_{S/B}}    
	\ar"1,1";"1,2"
	\hole
}
\end{equation} 
where we denote the restriction with $\phi$ itself, with a little abuse of notation. Moreover, a direct computation shows that the image of $\omega_B\otimes \cK_\partial$ via the morphism $\phi$ is contained $\cK_\partial\hookrightarrow f_*\omega_{S/B}$.

The Massey-product $\fm_{\sigma}(s_1,s_2)\in \Gamma(A, f_*\omega_{S/B}) $ of the pair of sections $s_1,$ $s_2$ of $\cK_\partial$ over a subset $A$ of $B$ is computed by \ref{Mor-MPRelativeFormsKernel} modulo the the $\cO(A)-$submodule $<s_1,s_2>$ of $\Gamma(A, f_*\omega_{S/B}).$

         	\subsection{Liftings of sections of $\cK_\partial$ with  Massey-trivial products.}
         	We assume that the rank $\rk\cK_\partial\geq 2$. We consider an open subset $A\subset B$ and a subspace of (local) sections  $W\subset\Gamma(A,\cK_\partial )$ such that $\dim_\bC W\geq 2.$ 
   	       	 \begin{definition} \label{Def-MTSubspSec}A subspace $W\subset\Gamma(A,\cK_\partial )$ is {\em Massey-trivial } if each pair of sections on $W$ is Massey-trivial (Definition \ref{Def-MPAJLoc}).
         	 \end{definition}

We can of course use the splitting given in \ref{SeS-KernelSplit<} to lift $W$ to $f_*\Omega^1_S$ and then we apply morphism \ref{Mor-MPRelativeFormsKernel} to compute the Massey-products of each pair of $W.$ This procedure defines sections of $f_*\omega_{S/B}$ that lie in $<s_1,s_2>_{\cO_B(A)},$ by definition of Massey-trivial pairs. We now prove that one can choose suitable liftings on $f_*\Omega^1_S$ with wedge zero (even if maybe different from those given by the splitting \eqref{SeS-KernelSplit<}).
\begin{proposition} \label{Prop-MTtoIsotropic} Let $A$ be an open contractible set of $B$ and $W\subset \Gamma(A,\cK_\partial)$ be a Massey-trivial subspace of sections of $\cK_\partial$ on $A.$ Assume that the evaluation map $W\otimes \cO_A\to {\cK_\partial}_{|A}$ defines an injective map of vector bundles. Then there exists a unique $\widetilde{W}\subset H^0(A,f_*\Omega^1_S)$ which lifts $W$ to $f_\ast\Omega^1_S$ and such that 
	the map \ref{Mor-MPRelativeForms} vanishes identically over $\widetilde{W},$ namely
	\begin{equation}\label{Mor-MPRlift}
	\xymatrix@!R{
		{\tilde{\phi}: \bigwedge^2\widetilde{W}\otimes {T_A}}                &  {{f_\ast\omega_{S/B}}_{|A,}}    
		\ar"1,1";"1,2"
		\hole
	}
	\end{equation} 
	is the null map.          	 		
\end{proposition}	
\begin{proof}  
	Let $\tau\in  T_A$ be a trivialization of $T_A$  and set $\beta$  its dual (that is, $\tau\cdot \beta=1$). 
	By composition of $W\otimes \cO_A\to {\cK_\partial}_{|A}$ with the splitting $\cK_\partial \to f_\ast \Omega^1_S,$ we obtain a lifting map $\rho : W\to \Gamma(A,  f_\ast \Omega^1_S).$
	Let $V=\rho(W)$ be its image. Let  $\{s_1,s_2,\dots,s_n\}$ be a basis of $W$ and set $v_i=\rho(s_i).$
	Since the Massey products are zero on any pairs of sections of $W,$ we have
	$$\tilde{\phi} (v_1\wedge v_i,\sigma)= f_is_1+g_is_i \quad \quad\mbox{and}\quad \quad \tilde{\phi} ( v_1 \wedge \sum _2^nv_i,\sigma)=f_0s_1+g_0 (\sum _2^ns_i),$$
	where the $f_i$ and the $g_i$ are holomorphic functions on $A.$ Assume $n=2$ and set
	$\tilde v_1=v_1-g_2\beta$ and
	$\tilde v_2= v_2-f_2\beta.$ Then $\tilde{\phi} (\tilde v_1\wedge  \tilde v_2,\sigma )=0$ and therefore 
	$\tilde v_1\wedge  \tilde v_2=0\in\Gamma(A, f_\ast \omega_S) .$ 
	We remark that  the unicity of the liftings $ \tilde v_1$ and $\tilde v_2$ follows at once.
	
	Assume by induction on $n>2$ that the proposition holds for $k<n.$ Consider the space $W'$
	generated by $\{s_1\dots s_{n-1}\}$ where we find liftings $\tilde v_i,$ for $i=1,\dots n-1,$ such that $\tilde v_i\wedge \tilde v_j=0.$ Moreover,  
	$$\tilde{\phi}(\tilde v_1\wedge v_n,\sigma)= es_1+fs_n \quad \quad\mbox{and}\quad \quad  \tilde{\phi} (\tilde v_1 \wedge( \sum _2^{n-1}\tilde v_i+ v_n),\sigma)=\tilde{\phi}  (\tilde v_1\wedge v_n,\sigma)=
	ls_1+m (\sum _2^ns_i).$$
	We get  $es_1+fs_n=ls_1+ m(\sum _2^ns_i)$ and $e=l$ and $f=0=m,$ since the sections are independent.
	Set $\tilde v_n=v_n-e\beta.$ Then $\tilde{\phi} ( \tilde v_1\wedge\tilde v_n,\sigma)=0$ and thus also
	$\tilde v_1\wedge\tilde v_n=0.$ In the same way we obtain $\tilde{v_i}\wedge \tilde{v_n}=0$. Unicity of the lifting follows immediately.

\end{proof}

\section{Flat sections and liftings to the sheaf of closed holomorphic forms on $S.$}\label{Sec-LiftingsOnU}
 
 In this section we relate the local system $\bU$ underlying the unitary flat bundle $\cU$ in the second Fujita decomposition of a fibration $f:S\to B$ with the subsheaf $\Omega^1_{S,d}\subset \Omega^1_S$ of the closed holomorphic $1$-forms on $S.$ We prove that it lifts to the direct image of the sheaf $\Omega^1_{S,d}$ so that it is described by closed holomorphic forms defined on tubular neighborhoods a  fibres of $f$. 
  
 \subsection{ Relative holomorphic de-Rham and a useful short exact sequence.}\label{SubSec-RelHoldeRham} 
 Given a fibration $f:S\to B$ there is a suitable short exact sequence which can be constructed by comparing the holomorphic de Rham sequences of the surface $S$ and the curve $B.$  

 Let us consider the holomorphic de-Rham sequence on $S$
 \begin{equation}\label{Ses-HoldeRahmS}
 \xymatrix@!R{
 	0  &    {\bC_S}    &    {\cO_S}    & {\Omega^1_{S,d}}   &   {0} ,    
 	\ar"1,1";"1,2"\ar"1,2";"1,3"\ar"1,3";"1,4"^\de\ar"1,4";"1,5"
 	\hole
 }	 
 \end{equation}
 where $\Omega^1_{S,d}= \ker (d:\Omega^1_S\to \Omega^2_S)$ is the sheaf holomorphic $1$-forms on $S$ that are $d$-closed. By pushing forward, we can relate
  it to the holomorphic de-Rham sequence on $B$
 \begin{equation}\label{Ses-HoldeRahmB}
 \xymatrix@!R{
 	0  &    {\bC_B}    &    {\cO_B}    & {\omega_B}   &   {0} ,  
 	\ar"1,1";"1,2"\ar"1,2";"1,3"\ar"1,3";"1,4"^\de\ar"1,4";"1,5"
 	\hole
 }	 
 \end{equation}
using the natural morphisms $\bC_B\to f_*\bC_S$ and $\cO_B\to f_*\cO_S$ induced by $f.$ Since they are both isomorphisms in this case, we obtain a diagram
 \begin{equation}\label{Dia-HolDeRhams}
 \xymatrix@!R{
 	{0}  & {f_*\bC_S}  & {f_*\cO_S}  & {f_*\Omega^1_{S,d}}  & {R^1f_*\bC_S}      & {R^1f_*\cO_S}  & \\
 	{0}  & {\bC_B}      & {\cO_B}  & {\omega_B}  & 0,         &                                     & 
 	\ar"1,1";"1,2"\ar"1,2";"1,3"\ar"1,3";"1,4"^{\de}\ar"1,4";"1,5"\ar"1,5";"1,6"
 	\ar"2,1";"2,2"\ar"2,2";"2,3"\ar"2,3";"2,4"^{\de}\ar"2,4";"2,5"
 	\ar@^{=}"2,2";"1,2" \ar@^{=}"2,3";"1,3" \ar@{^{(}->}"2,4";"1,4" 
 	\hole
 }
 \end{equation}
 which induces an injective morphism $\omega_B\hookrightarrow f_*\Omega^1_{S,d}$ together with the short exact sequence
 \begin{equation}\label{SeS-HolRelDeRham}
 \xymatrix@!R{
 	{0}  & {\omega_B}  & {f_*\Omega^1_{S,d}}  & {\cD}  & {0,} 
 	\ar"1,1";"1,2"\ar"1,2";"1,3"\ar"1,3";"1,4"\ar"1,4";"1,5"
 	\hole
 }
 \end{equation}
 where $\cD$ denotes the image of the morphism $f_*\Omega^1_{S,d} \to R^1f_*\bC_S.$ We will refer to this sequence as the {\em Relative holomorphic de-Rham sequence}. 
 We have the following properties
 \begin{lemma}\label{Lem-hatDInj} Let $f:S\to B$ be a semistable fibration of curves and let $\cD$ be the sheaf defined by \eqref{SeS-HolRelDeRham}. Then we have the injections of sheaves 
 	$(i)$ $i_{\mtin{$r$}}:\cD\hookrightarrow R^1f_\ast \bC$ and $(ii)$ $i_{\mtin{$\cD$}}:\cD\hookrightarrow \cK_\partial$.
 	\end{lemma}
 
 \begin{proof} $(i)$ This follows imediately by looking at the first row of diagram \eqref{Dia-HolDeRhams}, which is exact. Since by definition $\cD$ is exaclty the image of the map $r:f_\ast\Omega^1_{S,d}\to R^1f_\ast\bC$, the restriction $i_{\mtin{r}}$ of this to $\cD$ is an injection. $(ii)$ We consider the natural injection of sheaves  $ i_d:\Omega^1_{S,d}\hookrightarrow\Omega^1_S$ and we compare sequences (\ref{SeS-HolRelDeRham}) and (\ref{SeS-DeRhamOnFibf_*}). 
	We get the following diagram 
	\begin{equation}\label{Dia-LinkDiagramShort}
	\xymatrix@!R{
		0 &       {\omega_B}       &   {f_*\Omega^1_{S,d}}                &  {\cD}       &   {0}     &        \\
		0   & {f_*f^*\omega_B}    & {f_*\Omega^1_{S}}  & {f_*\Omega^1_{S/B}}  
		& {R^1f_*\cO_S\otimes \omega_B,} &       &       
		\ar"1,1";"1,2"\ar"1,2";"1,3"\ar"1,3";"1,4"\ar"1,4";"1,5"
		\ar"2,1";"2,2"\ar"2,2";"2,3"\ar"2,3";"2,4"\ar"2,4";"2,5"
		\ar@^{=}"1,2";"2,2"       \ar@{^{(}->}"1,3";"2,3"^{i_{\de}}   \ar@{^{(}->}"1,4";"2,4"^{i_{\mtin{$\cD$}}}   \ar"1,5";"2,5" 
		\ar"2,4";"2,5"^-{\partial}  
		\hole
	}
	\end{equation}
	that defines the injection $i_{\mtin{$\cD$}}:\cD\hookrightarrow \cK_{\partial}$ as claimed. 
	
	\end{proof}	
	\subsection{The local system $\bU$ and the sheaf $f_*\Omega^1_{S,d}$}\label{SubSec-UandTubForms}
	
		 We now describe the sheaf of flat sections of $\cU$ (that is the local system $\bU$) in terms of closed holomorphic forms on $S$ (that is sections of $f_*\Omega^1_{S,d}$).
		\begin{lemma}[Lifting lemma]\label{Lem-CharFlatSecN}
			Let $f:S\to B$ be a semistable fibration of curves and $\bU$ be the local system underlying the unitary summand $\cU$ of the second Fujita decomposition of $f.$ Then, there is a short exact sequence of sheaves
			Then there is short exact sequence 
			\begin{equation}\label{SeS-LiftingLemma}
			\xymatrix@!R{
				{0}  & {\omega_B}  & {f_*\Omega^1_{S,d}}  & {\bU}  & {0.} 
				\ar"1,1";"1,2"\ar"1,2";"1,3"\ar"1,3";"1,4"\ar"1,4";"1,5"\ar@/^-1.0pc/"1,4";"1,3"_<<<{\eta'}
				\hole
			}
			\end{equation}
			that is split.
			Moreover, the above sequence remains exact
			over each proper open subset $A$ of $B.$
		\end{lemma}
		\begin{proof}
				 We consider sequence \eqref{SeS-HolRelDeRham} and we prove that if $f$ is semistable, then the sheaf $\cD$ is isomorphic to $\bU$, so we get the stated short exact sequence. To do this we prove that $(i)$ $\cD$ is a local system whose stalk is isomorphic to a subspace $D\subset H^0(\omega_{F})$, for a general fibre $F$ of $f$; $(ii)$ $D$ is the larger subspace of $H^0(\omega_F)$ and so by definition (see subsection \ref{SubSec-Prel-LocSystOnFibr}, or directly \cite{CatDet_TheDirectImage_2014}) it is $\bU$. 
				 
				 $(i)$ We consider the morphisms  $k_{\mtin{$\cD$}}:\cD\hookrightarrow R^1f_*\bC_S$ of Diagram \eqref{Dia-HolDeRhams} and  $\alpha:R^1f_*\bC\to j_*j^* R^1f_*\bC $ of Lemma \ref{Lem-LocInvIso}, which is an isomorphism since  $f$ is semistable. By composition, we define 
%
the morphism $\alpha_{\mtin{$\cD$}}: \cD\to j_*j^*R^1f_*\bC_S$, which is still injective since $\alpha$ is an isomorphism. We consider the injection $j: B^0=B\setminus B_0\hookrightarrow B$ of the locus $B^0$ parametrizing the smooth fibres of $f$. The restriction $j^*R^1f_*\bC$ is a  local system whose stalk on a point $b\in B^0$ is isomorphic to the first cohomology group $H^1(F,\bC)$ of the corresponding fibre $F$. Then as a subsheaf of a local system, the rstriction $j^*\cD\subset j^*R^1f_\ast\bC$ must be a local subsystem of $j^*R^1f_*\bC$. Moreover, the stalk on $b\in B^0$ must be isomorphic to a vector subspace $D\subset H^{1,0}(F)=H^0(\omega_F)$ because $D$ is defined by taking the image through the map $f_\ast \Omega^1_{S,d}\to R^1f_\ast \bC$, which is given by restriction to the fibre so it must be of type $(1,0)$ (we note that by the exactness it is also exactly the kernel of the map $R^1f_\ast\bC\to R^1f_\ast \cO_S$, where $R^1f_\ast \cO_S$ is locally free with general fibre $H^{0,1}(F_b)$ while the map defines the projection $H^1(F_b,\bC)\to H^{0,1}(F_b)$).  But now we have that $j^\ast\cD$ is a local system whose general fibre $D\subset H^{1,0}(F)\subset H^1(F,\bC)$ and so by definition (see subsection \ref{SubSec-Prel-LocSystOnFibr} or directly \cite{CD:Answer_2017}) it must be a unitary local subsystem of $j^\ast\bU$. By applying the same argument given for $\bU$ we have that by pushing forwanrd $j_\ast j^\ast\cD$ it extends trivially to a local system on $B$ which is isomorphic to $\cD$. 

				$(ii)$ We prove that $j^\ast\cD$ is the larger local subsystem of $j^*R^1f_*\bC$ whose general fibre $D$ lies in $H^{1,0}(F)$.
				Let $\bD'$ be a local system whose  stalk is contained in $H^{1,0}(F)$ on a general fibre $F$. Then the map $\bD' \to R^1f_*\cO_S$ is zero and therefore $\bD'\subset \bD.$ 
				 
				 We have proven the existence of the exact sequence. To conclude that induced the long exact sequence in cohomology is still exact on a proper open subset $A$ of $B$ it is enough to note that $A$ is Stein and so $H^1(A,\omega_B)$ is zero.

				 We now prove that the sequence is split. We consider the section $\eta  \in H^0(\cK_\partial^\vee\otimes f_*\Omega^1_S)$ that splits the sequence \eqref{SeS-KernelSplit<}
				 (see Lemma \ref{Lem-SesSPLITKER}) and we prove that it induces also a splitting on sequence \eqref{SeS-LiftingLemma}. We rewrite diagram \eqref{Dia-LinkDiagramShort}
				 \begin{equation}\label{Dia-LinkDiagramShortU}
				 \xymatrix@!R{
				 	0 &       {\omega_B}       &   {f_*\Omega^1_{S,d}}                &  {\bU}       &   {0}     &        \\
				 	0   & {\omega_B}    & {f_*\Omega^1_{S}}  & {\cK_\partial}  
				 	& {0} &       &       
				 	\ar"1,1";"1,2"\ar"1,2";"1,3"\ar"1,3";"1,4"\ar"1,4";"1,5"
				 	\ar"2,1";"2,2"\ar"2,2";"2,3"\ar"2,3";"2,4"\ar"2,4";"2,5"
				 	\ar@^{=}"1,2";"2,2"       \ar@{^{(}->}"1,3";"2,3"^{i_{\de}}   \ar@{^{(}->}"1,4";"2,4"^{i_{\mtin{$\bU$}}}  
				 	\hole
				 }
				 \end{equation}
				 The kernel of the two sequences is the same. So the morphism $\eta':\bU\to f_*\Omega^1_S$ given by composition of $\eta$ with the injection $i_{\mtin{$\bU$}}:\bU\hookrightarrow \cK_\partial$ has image  $i_d:f_*\Omega^1_{S,d}\hookrightarrow f_\ast\Omega^1_S$. Thus $\eta'$ gives the desired splitting.
		\end{proof} 

\begin{remark} \label{Rem-LiftAreClosed} We notice that \begin{itemize}\item[$(i)$]
	any local lifting of a section of $\bU$ given by the lemma is a closed holomorphic forms since two liftings of the same section differ from the pullback of a holomorphic form on the curve $B$; \item[$(ii)$] any global section of $\bU$ contributes to the trivial summand $\cO^{\oplus q_f}\subset \cU$ of rank the relative irregularity $q_f$ of $f$ and we find the same description given in \cite{Gonz_PhdTs_2013}; \item[$(iii)$] while \cite[Corollary $7.2$]{ChenLuZu_OnTheOort_2016} states a lifting for the whole $\cU$ (and thus also of $\bU$) to $f_*\Omega^1_S$, we focus on the local system and we describe it more precisely. \end{itemize}
 	\end{remark}

 				\section{Massey-trivial families, flatness and the Castelnuovo-de Franchis theorem for fibred surfaces}\label{Sec-MTsubbundles}
 				In this section we restate a generalized version of the  Castelnuovo de Franchis theorem for fibred surfaces (see \cite{GonStopTor-On}) and we use it to relate the geometry of Massey-trivial  subspaces $W\subset \Gamma(A,\bU)$ of flat local sections of the unitary summand $\cU$ to the existence of a fibration from the surface into a smooth compact curve $\Sigma$ of genus greater than $2,$ after a base change. 
 				
 				We consider a fibration $f:S\to B$ of curves on a smooth curve $B$ which we do not assume to be compact (so the same holds for the surface $S$) and we assume that a general fibre $F$ has genus $g\geq 2$. Since $S$ is a complex surface we can still consider the sheaf  $\Omega^1_{S,d}$ of closed holomorphic $1-$forms on $S$ together with the wedge map 
			\begin{equation}
			\wedge: \bigwedge^2H^0(S,\Omega^1_{S,d})\to H^0(S,\omega_S) .
			\end{equation}
			\begin{definition}We say that a subspace $V$ of $H^0(S,\Omega^1_{S,d})$ is {\em isotropic} if the $\wedge-$map restricts to the null map on $\bigwedge^2 V$ and that it is moreover {\em maximal} if it is not properly contained in any larger isotropic space. 
				\end{definition}
				Then the theorem is the following.

 				\begin{theorem}[Castelnuovo de Franchis for fibred surfaces]\label{Thm-CdFL} Let $f:S\to B$ be a fibred surface on a smooth curve $B$ and let $V\subset H^0(S,\Omega^1_{S,d})$ be a maximal isotropic subspace of dimension $r\geq 2$ such that the restriction $V \to H^0(\omega_F)$ to a general fibre $F$ is injective. Then there is a non constant morphism $\varphi : S\to \Sigma$ from the surface $S$ to a smooth compact curve $\Sigma$ of genus $g_{\mtin{$\Sigma$}}\geq 2$ such that $\varphi^*H^0(\omega_{\Sigma})=V.$
 				\end{theorem}
 				\begin{proof} Let us first assume that $B$ is a complex disk. In this case the theorem has been proven in \cite{GonStopTor-On}, following the classical argument given in \cite[Proposition X.6]{Bea96} for a compact $S$. This assumption is essentially used to conclude that the foliation defined by pointwise proportional holomorphic forms is integrable, since all holomorphic forms on compact surfaces are closed. By releasing it, one can assume that holomorphic forms be closed (that is no longer automatic) and  check that this is sufficient to repeat the original argument. Let's briefly recall this point. Let $\omega_1, \omega_1,\ldots,\omega_k$ be a basis of $V$. Since $V$ is isoptropic, then we can find meromorphic functions $g_2,\ldots,g_k$ on $S$ such that $\omega_i = g_i\omega_1$. By using   $d\omega_i=0$, we differentiate these and we obtain $dg_i \wedge \omega_1=0$, for $i = 2,\ldots,k$. Furthermore, there is also a meromorphic function $g_1$ such that $\omega_1 = g_1 dg_2$ and $ dg_1 \wedge dg_2=0$. This means that the meromorphic differentials $dg_2, \ldots, dg_k$ are pointwise proportional to $\omega_1$ and also to $dg_1$, hence $dg_i \wedge dg_j = 0$ for any $i,j$. These functions define a meromorphic map $\psi: S \dra \bP^k$ as $\psi(p) = \left(1:g_1(p):g_2(p):\ldots:g_k(p)\right)$ and after eventually shrinking the disk in such a way that the  indeterminacy locus of $\psi$ is finite, one can take a resolution getting a holomorphic map whose image is a curve and conclude following the classical argument. When $B$ is not a disk, then we can always work locally and the result follows by a standard argument of analytic continuation.\end{proof} 
 				
 				We now consider a fibration $f:S\to B$ of curves on a projective curve $B$, the unitary flat summand $\cU$ of its second Fujita decomposition and we take a subspace $W\subset \Gamma(A, \bU)$ of sections on an open subset $A$. We recall that in section \ref{Sec-LiftingsOnU} we have proven that $\bU$ lifts to $f_*\Omega^1_{S,d}$ (Lemma \ref{SeS-LiftingLemma}). We get the following
 				\begin{proposition} \label{Prop-UMTtoIsotropic} Let $A$ be an open set of $B,$ let $i_A: A\hookrightarrow B$ be the inclusion and let $W\subset \Gamma(A,\bU)$ be a Massey-trivial subspace of sections of $\bU$ on $A.$ Then there exists a unique $\widetilde{W}\subset H^0(A,f_*\Omega^1_{S,d})$ which lifts $W$ to $f_\ast\Omega^1_{S,d}$ and such that    $\bigwedge ^2\widetilde{W}\to \Gamma(A, f_\ast \omega_S)$ is the  zero map. 					
 				\end{proposition}	
 				\begin{proof} We choose liftings of $W$ to $f_*\Omega^1_{S,d}$ by using Lemma \ref{SeS-LiftingLemma}. By noting that the evaluation map $W\otimes \cO_A\to i_A^*\cK_\partial$ is automatically an injective map of vector bundles, since $\Gamma(A,\bU)$ is the space of flat sections of $\cU$, we work pointwise applying Proposition \ref{Prop-MTtoIsotropic} to construct liftings with wedge zero. Then these are also closed since by Lemma \ref{SeS-LiftingLemma} they must differs from the others by sections of $\omega_B.$ 
 					\end{proof}
 					We showed that Massey-trivial subspaces $W$ of sections on $\bU$ correspond to isotropic subspaces of sections on $f_\ast\Omega^1_{S,d}.$ Now we relate it to the monodromy of the flat bundle generated by $W$. Let $H<\pi_1(B,b)$ be the kernel of the monodromy representation of $\bU$ (see section \ref{SubSec-Prel-LocSystOnFibr}) and set $H_{\mtin{W}}$ be the subgroup of $H$ which acts trivially on $W,$ i.e. 
 					\begin{equation}
 					H_{\mtin{W}}=\{g\in H\, |  \, g\cdot w= w, \forall w\in W\}.
 					\end{equation}
 					
 					We obtain the following

 					\begin{theorem}\label{Prop-MPpencils2} Let $f:S\to B$ be a semistable fibration of curves of genus $g\geq 2$ on a projective curve $B$. We consider a subspace $W\subset \Gamma(A,\bU)$ of local flat sections on $A$ that is a maximal Massey-trivial subspace and a subgroup $K$ of $H_{\mtin{W}}$. Then the fibration  $f_{\mtin{K}}:S_{\mtin{K}}\to B_{\mtin{K}}$ given by the base change with $u_{\mtin{K}}: B_{\mtin{K}}\to B$, which is the morphism classified by $K$, has an irrational pencil $h_{\mtin{K}}:S_{\mtin{K}}\to \Sigma$ on a smooth compact curve $\Sigma$ such that $W\simeq h_{\mtin{K}}^*H^0(\omega_{\Sigma}).$ 
 						
 					\end{theorem}
 					
 					\begin{proof}
 						Let $u_{\mtin{W}}:B_{\mtin{W}}\to B$ be the \'{e}tale covering classified by  $H_{\mtin{W}}< H.$ By construction, the pull back of $W$ extends to a subspace $\widehat{W}$ of global sections $\Gamma(B_{\mtin{W}},\bU_{\mtin{W}}),$ where  $\bU_{\mtin{W}}$ is the unitary summand of the fibration $f_{\mtin{W}}:S_{\mtin{W}}\to B_{\mtin{W}}$ defined by the base change. By applying Proposition \ref{Prop-UMTtoIsotropic}, we get the proof for $\widehat{W}$. Then we can repeat the argument for any subgroup $K$ of $H_{\mtin{W}}$.
 						
 					\end{proof}
 					
%
 					Now let us note that a subspace $W\subset \Gamma(A, \bU)$ of local sections on a contractible $A$ is not necessarily invariant under the monodromy action of $\bU$, so we give the following (see also section \ref{SubSec-Prel-LocSystOnFibr}).
 					
 					\begin{definition}\label{Def-MTsubbundleB} 
 						We say that a local subsystem $\bM$ of $\bU$ is {\em Massey-trivial} if its stalk $M$ is isomorphic to a Massey-trivial subspace of  $\Gamma(A,\bU)$ of sections an open contractible $A\subset B$; we say that it is {\em Massey-trivial generated} if the stalk $M$ is generated by a Massey-trivial subspace of $\Gamma(A,\bU).$  
 					\end{definition}  
 					
 					\begin{remark}\label{rem-MTaction} By using a standard argument of analytic continuation, we have that the property above is stable under the monodromy action on subspaces $W \subset \Gamma(A,\bU)$ of local flat sections on $\bU$. 
 						\end{remark} 				


	\section{Proof of the main theorems} \label{Sec-ProofMainTheorems}

In this section we prove our main results (Theorems \ref{Thm-MainG} and \ref{Thm-MainSbis} and Corollary \ref{cor:main}). 
Let $f:S\to B$ be a fibration of curves of genus $g\geq 2$ on a smooth projective curve $B$, let $b$ be a general point of $B$ and let $F$ be the corresponding fibre (here general means that $F$ is smooth and is non special with respect to the Massey-trivial vanishing property, e.g. hyperelliptic unless it is not the case for all the fibres). We recall that $\cU$ is the unitary flat summand of the second Fujita decomposition of $f$, $\bU$ is the underlying local system (i.e $\cU=\bU\otimes \cO_B$) and $\rho:\pi_1(B,b)\to U(	\rank{$U$}, \bC)$ is its monodromy representation, where $\rank{$U$}$ is the rank of $\cU$, $H=\ker \rho$ is the kernel and $G=\pi_1(B,b)/H$ is the quotient that we identify with the monodromy group as usual. We consider a flat subbundle $\cM\subset \cU$ that is generated by a (maximal) Massey-trivial subspace. Then in the same way, $\bM\subset \bU$ is the underlying local subsystem, $\rho_{\tiny{M}}:\pi_1(B,b)\to U(\rank{$M$},\bC)$, $\rho_{\tiny{M}}\subset \rho$, is its monodromy subrepresentation, $H_{\mtin{$M$}}$ is the kernel and $G_{\mtin{M}}=\pi_1(B,b)/H_{\mtin{M}}$ is the monodromy group.
 We recall that $\cM$ is generated by a maximal Massey-trivial vector space if there is a subspace $W\subset \Gamma(A,\bU)$ of local sections such that $W$ is Massey-trivial, generates $\cM$, i.e. $\bM$ has stalk $M=G\cdot W$ and rank $\rank{$M$}=\dim M$, and $W$ is not properly contained in any other subset of $\Gamma(A,\bU)$ that is both generating and Massey-trivial.  

\subsection*{Proof of Theorem \ref{Thm-MainG}} We assume that $f$ is semistable and we prove that in this case the monodromy acts faithfully on a finite set of morphisms from a general fibre $F$ to a fixed smooth projective curve $\Sigma$ of genus $g_{\mtin{$\Sigma$}}\geq 2$. We have two steps: in step $1$ we define our family of morphisms while in step $2$ we define our faithful action.

{\bf Step $1$: construction of a set of morphisms of curves.}
We consider the Galois covering map $u_{\mbox{{\tiny $\cM$}}}:B_{\mbox{{\tiny  $\cM$}}}\to B$  classified by the normal subgroup $H_{\mtin{$\cM$}}$ of $\pi_1(B,b)$. This is defined by the action of the monodromy group $G_{\mbox{{\tiny  $\cM$}}}=\pi_1(B,b)/ H_{\mbox{{\tiny  $\cM$}}}$ of $\cM$ and so it trivializes the monodromy of $\cM$. We consider the \'{e}tale base change  
\begin{equation}\label{Dia-BaseChangeM}
\xymatrix@!R{
	{S_{\mtin{$\cM$}}:=S\times_BB_{\mtin{$\cM$}}} \,        &     {S}              &    \\
	{B_{\mtin{$\cM$}}}            \,        & {B}         &    
	\ar"1,1";"2,1"^{f_{\mtin{$\cM$}}}   \ar "1,1";"1,2"^>>>>>>{\varphi_{\mtin{$\cM$}}}
	\ar"1,2";"2,2"^{f} \ar "2,1";"2,2"^{u_{\mtin{$\cM$}}}
	\hole
}
\end{equation}
getting a fibration $f_{\mtin{$\cM$}}:S_{\mtin{$\cM$}} \to B_{\mtin{$\cM$}}$. Then $S_{\mtin{$\cM$}}$ is a smooth surface given by the fibred product $S\times_{B}B_{\mtin{$\cM$}}$ (that is non compact  when the monodromy is infinite) and $\varphi_{\mtin{$\cM$}}:S_{\mtin{$\cM$}}\to S$ is a fiber preserving \'{e}tale Galois covering. By construction, the action of $G_{\mtin{$\cM$}}\times S_{\mtin{$\cM$}}\to S_{\mtin{$\cM$}}$ sends a point $(p,b)$ to the point $(p, gb'),$ for $g\in G_{\mtin{$\cM$}}$, under the action of the automorphism $b'\mapsto gb'$ of $B_{\mtin{$\cM$}}$ defined by the action of $G_{\mtin{$\cM$}}$ on $B_{\mtin{$\cM$}}.$ We set $g(p,b')=(p, gb')$ and note that $g: S_{\mtin{$\cM$}}\to S_{\mtin{$\cM$}}$ defines an automorphism of $S_{\mtin{$\cM$}}$ compatible with the fibration $f_{\mtin{$\cM$}}$ (i.e. fiber preserving). 

After identifying the maximal Massey-trivial subset $W\subset \Gamma(A,\bU)$ generating $\cM$ with its image inside $H^0(\omega_{F_b})$ under the evaluation map (Remark \ref{garibaldi}), a direct computation allows to describe $H_{\mtin{$\cM$}}$ as 
\begin{equation}\label{Eq-desKer}
H_{\mbox{{\tiny  $\cM$}}}=\{g\in \pi_1(B,b)| \, gg'w=g'w, \, \forall w\in W, \, \forall g'\in G_{\mtin{$\cM$}}\}.
\end{equation}
We can apply Proposition \ref{Prop-MPpencils2} (since $f$ is semistable) using  $H_{\mtin{$\cM$}}$ as a subgroup $K$ of $H_{\mtin{W}}$ and getting  a non-constant map $h:S_{\mtin{  $\cM$}}\to \Sigma$ on a smooth projective curve $\Sigma$ such that $u_{\mtin{$\cM$}}^*W\simeq h^*H^0(\omega_{\Sigma}).$ Under the action of the group $G$, the map above constructs a set $$\mathscr{H}:=\{h_i:S_{\mtin{  $\cM$}}\to \Sigma\, | \, i\in \cI\}$$ of non constant morphisms. By restriction to $F$, we obtain a set  $$\mathscr{K}:=\{k_i:F\to \Sigma\, | \, i\in \cI\}$$ of morphisms between curves of genus greater than $2$ parametrized by a set $\cI$ of indices. More precisely, we fix a point $b_0$ in the fiber of $b$ with respect to $u_{\mtin{$\cM$}}$ and we denote with $F_0$ the fibre of $f_{\mtin{$\cM$}}$ over $b_0$ (that is naturally isomorphic to $F$). For any $g\in G_{\mtin{$\cM$}},$ we perform the automorphism $g:S_{\mtin{$\cM$}}\to S_{\mtin{$\cM$}}$ to define $h_g$ and $k_g$ by composition as in the following
\begin{equation}\label{Dia-PencMor}
\xymatrix@!R{
	{F_0}  & {S_{\mtin{$\cM$}}} \,        &     {S_{\mtin{$\cM$}}}       &      \\
	& {}                               &     {\Sigma}
	\ar@{^{(}->}"1,1";"1,2"^{i}   \ar "1,2";"1,3"^{g} 
	\ar"1,3";"2,3"^{h} 
	\ar"1,1";"2,3"_{k_g} \ar "1,2";"2,3"^{h_g}
	\hole
} 						
\end{equation}
In general $\mathscr{K}$ is not parametrized by the whole $G_{\mtin{$\cM$}}$ since two elements should define the same map. For example this happens when there is $g\in G_{\mtin{$\cM$}}$ that fixes not only $W$ but every element of $W$ under the monodromy action. In this case $k_g=k$. Anyway by construction we can fix an injection $\cI\hookrightarrow G_{\mtin{$\cM$}}$ sending $i\in \cI$ to $g_i\in G_{\mtin{$\cM$}}$.

{\bf Step 2: construction of a faithfull action.}
We consider the actions on $\mathscr{H}$ and $\mathscr{K}$
\begin{equation}\label{Act-IrrPenc}
\xymatrix@!R{
	{G_{\mtin{$\cM$}}\times \mathscr{H}}  & {\mathscr{H}} ,\,        &     {(g_1,h_{i})}       &  {g_1\cdot h_{i}:=h_{j}, \, \,j\in \cI \,:\, h_j=h_{g_1g_i}}  							
	\ar"1,1";"1,2"
	\ar@{|->} "1,3";"1,4" 
	\hole
} 	
\end{equation}
\begin{equation}\label{Act-Mor}
\xymatrix@!R{
	{G_{\mtin{$\cM$}}\times \mathscr{K}}  & {\mathscr{K}} ,\,        &     {(g_1,k_{i})}       &  {g_1\cdot k_{i}:=(g_1\cdot h_{i})\circ i=k_{j}, \, \,j\in \cI \,:\, k_j=k_{g_1g_i}}  							
	\ar"1,1";"1,2"
	\ar@{|->} "1,3";"1,4" 
	\hole
} 	
\end{equation}
naturally defined by the action of $G_{\mtin{$\cM$}}$ on $S_{\mtin{$\cM$}}$. From now on we will write $g_1g_i$ instead of $j$ with a little abuse of notation. We set  $\Bij(\mathscr{K})$ the group of the bijections on $\mathscr{K}$ and we consider the homomorphism 
\begin{equation}\label{Map-AutK}
\xymatrix@!R{
	{\Psi_{\mtin{$\cM$}}\,\colon\,G_{\mtin{$\cM$}}}  & {\Bij(\mathscr{K})} ,\,        &     {g_1}       &  {g_1\cdot\,:\, \mathscr{K}\to \mathscr{K}} 
	\ar"1,1";"1,2"
	\ar@{|->} "1,3";"1,4" 
	\hole
} 	
\end{equation}
naturally induced by the action. 
We want to prove that $G_{\mtin{$\cM$}}$ acts faithfully on $\mathscr{K}$, which is equivalent to prove that $\Psi_{\mtin{$\cM$}}$ is injective. To do this, we need the following formula.

\begin{lemma}\label{Lem-formulaFormsMonodromy} Let $e\in G_{\mtin{$\cM$}}$ be the neutral element and take $\alpha\in H^0(\omega_{\Sigma}).$ Then for any $g\in G_{\mtin{$\cM$}},$
	\begin{equation}\label{Formula-PullBackMonod}
	{k^*_{g}(\alpha) = g^{-1} k^*_e(\alpha)}, 
	\end{equation}
	where $g^{-1}$ acts on $k^*_e(\alpha)\in W$ under the monodromy action $\rho_{\mtin{$\cM$}}$ of $\cM.$ 
\end{lemma}
\begin{proof}
	Let $A_0$ be an open contractible subset of $B_{\mtin{$\cM$}}$ around $b_0$ such that $u_{\mtin{$\cM$}}(A_0)\subset A.$ 
	Let $\widetilde{W}\subset \Gamma(A_0,{f_{\mtin{$\cM$}}}_\ast \Omega^1_{S_{\mtin{$\cM$}},d})$ be the unique lifting of $W$ provided by Proposition \ref{Prop-UMTtoIsotropic}. By construction, we can naturally lift the monodromy action from $W$ to $\widetilde{W}$ and then for $\eta\in \widetilde{W}$ we have $g\eta=\eta$, for any $g\in H_{\mtin{$\cM$}}.$ Namely, $\widetilde{W}$ extends to a subspace of global forms $H^0(S_{\mtin{$\cM$}}, \Omega^1_{S_{\mtin{$\cM$}},d})$ and we identify them. Let $\eta=h_e^*\alpha\in \widetilde{W}$ and let $w=\eta_{|F}.$ Then we have that $k^*_{g}(\alpha)=(g^*h_e^*\alpha)_{|F}=(g^*\eta)_{|F}=\eta_{|F_{g^{-1}b}}$ and $g^{-1}k^*_{e}(\alpha)=g^{-1}(\eta_{|F})=g^{-1}w=\eta_{|F_{g^{-1}b}}$ are the same.
\end{proof}

We now apply the previous formula to prove that the action is faithful, showing the following
\begin{proposition}\label{Lem-PsiInj} The map $\Psi_{\mtin{$\cM$}}$ is injective.
\end{proposition}
\begin{proof}
	We consider the neutral element $e$ of $G_{\mtin{$\cM$}}$ and $g_1 \in G_{\mtin{$\cM$}}$ such that $g_1\neq e$ and we prove that $\Psi_{\mtin{$\cM$}}(g_1)\neq \Psi_{\mtin{$\cM$}}(e).$ By definition we have to show that there is $g_2\in G_{\mtin{$\cM$}}$ such that $g_1\cdot k_{g_2}\neq  e\cdot k_{g_2}.$ Let $\mathscr{M}(F,\Sigma)$ be the set of morphisms between $F$ and $\Sigma$ and let  $\Hom(H^{1,0}(\Sigma),H^{1,0}(F))$ be the set of homomorphisms between the spaces of their holomorphic forms. Since the pullback functor $\mathscr{M}(F,\Sigma)\to \Hom(H^{1,0}(\Sigma),H^{1,0}(F))$ is injective (see e.g. \cite{Martens_Obervations_1988} and also \cite{AlzatiPirola_Some_1991}), we equivalently prove that $k_{g_1g_2}^*\neq k_{g_2}^*$. To do this, we apply the description of $H_{\mtin{$\cM$}}$ given in \ref{Eq-desKer}. Since $g_1\neq e$, then $g_1\notin H_{\mtin{$\cM$}}$ and so we get $w\in W\subset H^0(\omega_F)$ and $g_2\in G_{\mtin{$\cM$}}$ such that $g_1g_2w\neq g_2w.$ Let $\alpha \in H^0(\omega_\Sigma)$ be such that $k_e^*\alpha=w.$ We apply Formula \ref{Formula-PullBackMonod} to $g=g_1g_2,$ obtaining $g_1g_2w=(g_1g_2)^{-1} k_e^*(\alpha)=k^*_{g_1g_2}w.$ Then we also apply the same formula to $g=g_2$, obtaining $g_2w=g_2^{-1} k^*_e(\alpha)=k^*_{g_2}(\alpha).$ Since by assumption $g_1g_2w\neq g_2 w$, we conclude that  $k^*_{g_1g_2}\neq k^*_{g_2}.$
\end{proof}


\subsection*{Proof of Theorem \ref{Thm-MainSbis}} We work under the assumptions and notations of Theorem \ref{Thm-MainG}. 

 {\em Proof of  $(i)$ }We prove that the monodromy group $\cG_{\mtin{$\cM$}}$ is finite. To do this we apply a classical de-Franchis' theorem (see e.g.  \cite{Martens_Obervations_1988}) to the set of $\mathscr{M}(F,\Sigma)$ of morphisms between $F$ and $\Sigma$.  The theorem states that the set of non constant morphisms between two curves of genus greater than $2$ is finite. By Theorem \ref{Thm-MainG} we conclude the same for the monodromy group of $\cM$.

{\em Proof of $(ii)$.} We consider the base change $u_{\mbox{{\tiny M}}}:B_{\mbox{{\tiny M}}}\to B$ constructed in the proof of Theorem \ref{Thm-MainG}, which trivializes the monodromy group of $\cM$. Since the fiber $M$ of $\cM$ is generated by $W$ (i.e. $M=G\cdot W=\sum_{i\in \cI}g_i\cdot W$), we compute the pullback under $u_{\mbox{{\tiny M}}}$ using the family of morphisms parametrized by $\cI,$ for each $i\in \cI$ and getting $k_{g_i}^*H^0(\omega_\Sigma)=i^*h_{g_i}^*H^0(\omega_\Sigma)=g_iW$. Then $\cM$ has general fibre isomorphic to $\sum_{a\in G} k_{g_i}^*H^0(\omega_{\Sigma})\subset H^0(\omega_F)$.

\subsection*{Proof of Corollary \ref{cor:main}}
					 We consider $f:S\to B$ and $\cM\subset \cU$ and we apply the semistable reduction theorem (see e.g. \cite{CD:Answer_2017}) getting a finite morphism of curves $u:B'\to B$ that produces by base change a semistable fibration $f':S'\to B'$. Let $\cU'$ be the unitary flat summand of the second Fujita decomposition of the semistable fibration $f'$. By Lemma \ref{Lem-UunderBC}, we get an inclusion $u^\ast \cU \subset \cU'$ that is compatible with the underlying structure of local systems.
					 
					 We consider the pullback $\cM'=u^\ast\cM\subset \cU'$, which is of course still defined by a local subsystem of $\cU'$. We prove that the monodromy group $G_{\mtin{$\cM'$}}$ of $\cM'$ is finite. We can assume that $\cM'$ is generated by a maximal Massey-trivial subspace. 
Then by Theorem \ref{Thm-MainSbis}, $\cM'$ has finite monodromy and so by proposition \ref{Prop-SubLocSystMon}, we conclude that also $\cM$ has finite monodromy.
 					 			\section{Applications}\label{Sec-Applications}
 					 			\subsection{Monodromy and the Griffiths infinitesimal invariant of the canonical normal function}\label{SubSec-NonVanCriteria}
 					 			
 					 			In this section we relate the monodromy of the unitary flat summand $\cU$ of a fibration of curves to the Griffiths ifnfinitesimal invariant of the canonical normal function of a general fibre, when it is generated by a Massey trivial subspace. We use \cite[Theorem 3.2.1]{C-P_TheGriffiths_1995} together with Theorem \ref{Thm-MainSbis}.  
 					 			
 					 				        Let $f:S\to B$ be a fibration of curves of genus $g\geq 2$ on a smooth projective curve $B$, let $F$ be a general fibre and let $J(F)$ be its Jacobian. 
 					 				        We consider the weight-$(2g-3)$ geometric HS $(H_\bZ=H^{2g-3}(J(F),\bC), F^{p}H^{2g-3}(J(F),\bZ))$  of $J(F)$, which is given by
 					 				        $H^{p,q}=\bigwedge^p H^{1,0}\otimes \bigwedge H^{0,1}, \quad H^{1,0}\simeq H^0(\omega_F)$
 					 				        since $J(F)$ is a principally polarized abelian variety of dimension $g$.
 					 				        
 					 				        The $(g-1)-$Griffiths intermediate Jacobian of $J(F)$ is defined by the  weight-$(2g-3)$ geometric HS of $J(F)$ as
 					 				        \begin{equation*}
 					 				        J^{g-1}(J(F))=H^{2g-3}(J(F),\bC)/(F^{g-1}H^{2g-3}(J(F),\bC)\oplus H^{2g-3}_{\bZ})\simeq F^2H^3(J(F),\bC)^*/H_3(J(C),\bZ),
 					 				        \end{equation*} 
 					 				        where the isomorphism is given by the Poincar\'{e} duality.
 					 				        Moreover, if $[\Theta] \in H^2(J(F), \bZ)$ is the principal polarization on $J(F)$ given by the $\Theta$-divisor, the associated Lefschetz operator induces a decomposition 
 					 				     $   H^3(J(F),\bC)=P^3(J(F),\bC)\oplus [\Theta]\cdot H^1(J(C),\bC,$
 					 				        where $P^3(J(F),\bC)$ is the primitive cohomology. This defines the intermediate primitive cohomology 
 					 				        \begin{equation}
 					 				        P(F)=F^2P^3(J(F),\bC)^*/H_3(J(C),\bZ)_{prim}, 
 					 				        \end{equation} 
 					 				        where $H_3(J(F),\bZ)_{prim}$ is the image of $H_3(J(F),\bZ)$ in $F^2P^3(J(F),\bC)^*.$
 					 				        
 					 				        Let $Z^{g-1}(J(F))_{\mtin{hom}}$  be the group of one dimensional algebraic cycles homologically equivalent to zero in $J^{g-1}(F)$. The Abel Jacobi map
 					 				        $\phi\,: \, Z^{g-1}(J(F))_{\mtin{hom}}\to J^{g-1}(F)$
 					 				        is the map given by integration of cycles, depending on a base point $b\in B.$ The composition of the map above with the projection $J^{g-1}(F)\to P(F)$  is independent of the choice of the base point.
 					 				        
 					 				        The {\em Ceresa cycle} of the general fibre $F$ is defined as the one cycle $[F-F^{-}]\in Z^{g-1}(J(F))_{\mtin{hom}}$ given by the image of $F-F^{-}$ in $J(F)$ via the Abel Jacobi map.

 					 				        On the locus $j:B^0\hookrightarrow B$ parametrizing the smooth fibres of $f$ we can construct a fibration $\cP(f):\cP\to B^0$ of the primitive intermediate Jacobians associated to the Hodge structure 
 					 				        $H^{2g-3}(J(F),\bC)$ of the Jacobian $J(F)$ of the fibres $F$ of $f.$ 
 					 				        
 					 				        The canonical normal function $\nu: B^0\to \cP$ is the section associating to each $b\in B^0$ the image of the Ceresa cycle $[F_b-F_b^{-}]\in Z^{g-1}(J(F_b))_{\mtin{hom}}$ in $\cP$ through the higher Abel-Jacobi map composed with the projection over $\cP$. 
 					 				       
 					 				       The associated { Griffiths infinitesimal invariant }$\delta(\nu)$ (see e.g. \cite{C-P_TheGriffiths_1995} for the explicit definition) induces pointwise over a point $b\in B^0$ a linear map $\ker (\gamma) \to \bC $ from the kenrle of the natural map $\gamma : T_{B,b}\otimes P^{2,1}J(F_b)\to P^{1,2}J(F_b)$ defined by the IVHS on the primitive cohomological groups $P^{1,2}J(F_b)$ and $P^{1,2}J(F_b)$ in $H^3(J(F_b),\bC).$ 
 					 				       The formula provided in \cite[Theorem 3.2.1]{C-P_TheGriffiths_1995} depends on the description above, on the polarization  $Q(-,=)=\int_{F_b}-\wedge=$ of the HS $(H_\bZ=H^1(F_b,\bZ), H^{1,0}=H^0(\omega_{F_b})),$  given by the intersection form, and on the Kodaira-Spencer class $\xi_b\in H^1(T_{F_b})$  of the fibre $F_b.$ 
 					 				       
 					 				       \begin{lemma}\cite[Theorem 3.2.1]{C-P_TheGriffiths_1995}\label{Lem-GrifFor} Let $\omega_1,\omega_2,\sigma\in H^0(\omega_{F_b})$ be such that $\xi_b\cdot\omega_1=\xi_b\cdot\omega_2=0$ and $Q(\omega_1,\bar \sigma)=Q(\omega_2,\bar \sigma)=0.$ Then $\omega_1\wedge\omega_2\wedge \bar \sigma$ lies, up to a natural isomorphism, in $\ker \gamma$ and we have
 					 				       	\begin{equation} \label{For-CompGII}
 					 				       	\delta(\nu)(\xi_b\otimes \omega_1\wedge\omega_2\wedge \bar \sigma)= -2Q (\fm_{\xi_b}(\omega_1,\omega_2), \bar \sigma).
 					 				       	\end{equation} 
 					 				       \end{lemma}
 					 				        %
 					 				        
 					 				        Putting the formula together with Theorem \ref{Thm-MainG}, we have the following 
 					 				        \begin{corollary}\label{Cor-InfMon&NonVanGriff}
 					 				        	Let $f:S\to B$ be a fibration of curves of  genus $g\geq 2$ and let $\cU$ be the unitary summand in the second Fujita decomposition of $f.$ If the monodromy of $\cU$ is infinite, then the Griffiths infinitesimal invariant of the canonical normal function  $\nu: B^0\to \cP$ is non zero on a general point $b\in B^0$ and so the normal function $\nu$  is non torsion.
 					 				        \end{corollary}
 					 				        \begin{proof}
 					 				        	We apply formula \ref{For-CompGII} to sections of $j^\ast\cU\subset j^\ast\cK_\partial.$  Since the monodromy of $\cU$ is infinite, then by Theorem \ref{Thm-MainSbis} $\cU$ is not Massey-trivial generated and we can find a pair $(\omega_1,\omega_2)\subset H^0(\omega_{F_b})$ of independent elements such that  $\fm_{\xi_b}(\omega_1,\omega_2)\neq0,$ where $\xi_b$ is the Kodaira spencer class of $F_b.$ Applying formula \ref{For-CompGII} to $\omega_1,\omega_2\in H^0(\omega_{F_b})$ (which are such that $\xi_b\cdot\omega_1=\xi_b\cdot\omega_2=0$) and $\sigma=\fm_{\xi_b}(\omega_1,\omega_2)\in H^0(\omega_{F_b})$ 
 					 				        	we get 
 					 				        	\begin{equation}
 					 				        	\delta(\nu)(\xi_b\otimes \omega_1\wedge\omega_2\wedge \bar \sigma)= -2Q (\fm_{\xi_b}(\omega_1,\omega_2), \bar \fm_{\xi_b}(\omega_1,\omega_2))\neq 0
 					 				        	\end{equation} 
 					 				        	So the Griffiths infinitesimal invariant in not zero. To conclude, the normal function is consequently non torsion as proven in \cite{Grif_InfinitesimalVariationsIII_1983}, \cite{Green_1989},\cite{Vois_UneRemark_1988}.
 					 				        \end{proof}
 					 				        
 					 				        In particular, the previous result applies to the examples provided in \cite{CatDet_TheDirectImage_2014}, in \cite{CD:Answer_2017} and also in \cite{CatDet_Vector_2016}, concerning the construction of fibrations where the monodromy of $\cU$ is infinite. 
 					 				        \begin{corollary}\label{Cor-CataDetNFnotTor}
 					 				        	Let $f:S\to B$ be a fibration as in \cite{CatDet_TheDirectImage_2014}, \cite{CD:Answer_2017} and \cite{CatDet_Vector_2016} with $\cU$ of not finite monodromy. Then the canonical normal function is non torsion.
 					 				        \end{corollary}
 					 				
%
 					 			\subsection{Semiampleness and Massey-trivial products}\label{SubSec-SemiamplenessCriteria}
 					 			In this section we state a criterion for the semiampleness of $f_*\omega_{S/B}$ depending on Massey-trivial generated flat bundles. This is a corollary of Theorem \ref{rem-MTaction} together with a characterization of semiampleness for unitary flat bundles ( see e.g \cite[Theorem 2.5]{CD:Answer_2017}). In particular we show that the criterion applies to hyperelliptic fibrations.
 					 			
 					 			We recall that by the second Fujita decomposition $f_\ast\omega_{S/B} = \cU\oplus \cA$ of $f$ it follows that semiampleness depends only on the unitary flat bundle $\cU$. By \cite[Theorem 2.5]{CD:Answer_2017} a unitary flat bundle  is semiample if and only if it has finite monodromy. 
 					 			
 					 			Putting together with theorem \ref{rem-MTaction}, we have the following
 					 			
 					 			\begin{corollary}\label{Cor-SemiamplenessCriteria} Let $f:S\to B$ be a semistable fibration of curves of genus $g\geq 2$. If the unitary bundle $\cU$ of its second Fujita decomposition is Massey-trivial generated, then $f_*\omega_{S/B}$ is semiample.
 					 			\end{corollary}
 					 			We note that semiapleness is not a numerical property and the corollary gives a local condition (Massey-triviality) that can be checked to get a global one.
 					 			This for example works on hyperelliptic fibrations, namely fibrations whose general fibre is a hyperelliptic curve, giving another proof of a result stated in \cite{LuZuo_OnTheSlope_2017}.
	 					 			
 					 			 \begin{corollary} The sheaf $f_\ast \omega_{S/B}$ of a fibration $f:S\to B$ of hyperelliptic curves of genus $g\geq 2$ on a smooth projective curve $B$ is semiample.
 					 			 \end{corollary}
 					 			 \begin{proof} By using the second Fujita decomposition of $f$, we just have to prove that its unitary flat summand $\cU$ is semiample. To do this we prove that $\cU$ is Massey-trivial and then we conclude by corollary \ref{Cor-SemiamplenessCriteria}. We consider  two independent vectors $s_1,s_2\in U\subset H^0(\omega_F)$ and we compute the Massey product $m_{\xi}(s_1,s_2 )$ of the pair $(s_1,s_2).$ By formula \ref{Mor-Mp/Aj}, $m_{\xi}(s_1,s_2 )$ is antisymmetric with respect to $s_1,s_2.$ But now the hyperelliptic involution $\sigma: F\to F$ acts on $H^0(\omega_F)$ by pullback  $\sigma^*:H^0(\omega_F)\to H^0(\omega_F)$ as the $-1$ multiplication map so we get 
 					 			 	 $\sigma^*m_{\xi}(s_1,s_2 )=-m_{\xi}(s_1,s_2 )$ and $\sigma^*m_{\xi}(s_1,s_2 )=m_{\xi}(-s_1,-s_2 )=m_{\xi}(s_1,s_2 )$. Then it must be zero since both symmetric and antisymmetric.     
 					 			 \end{proof}
 					 		

%
%

\end{document}